\documentclass[reqno,final]{amsart}
\usepackage{natbib}  
\usepackage{fancyhdr} 
\usepackage{color} 
\usepackage{hyperref} 
\usepackage{graphicx} 
\usepackage[utf8]{inputenc}
\usepackage{enumerate}
\usepackage{amsmath, amssymb, amsthm,comment}
\usepackage{tikz}
\usepackage{dsfont}

\usepackage{pstricks}
\usepackage{amssymb}

\definecolor{aleacolor}{rgb}{0.16,0.59,0.78}

\hypersetup{
breaklinks,
colorlinks=true,
linkcolor=aleacolor,
urlcolor=aleacolor,
citecolor=aleacolor}

\renewcommand{\headheight}{11pt}
\def\emptyset{\varnothing} 

\def\d{{\rm d}} 
\def\e{\varepsilon}

\def\O{\Omega}

\def\L{\Lambda}

\newfam\Bbbfam 
\font\tenBbb=msbm10 
\font\sevenBbb=msbm7 
\font\fiveBbb=msbm5 
\textfont\Bbbfam=\tenBbb 
\scriptfont\Bbbfam=\sevenBbb 
\scriptscriptfont\Bbbfam=\fiveBbb

\newcommand{\R}     {\mathbb{R}} 
\newcommand{\Z}     {\mathbb{Z}} 
\newcommand{\N}     {\mathbb{N}} 
\renewcommand{\P}   {\mathbb{P}} 
 
\newcommand{\E}     {\mathbb{E}} 
\newcommand{\Q}     {\mathbb{Q}}

\newcommand{\smfrac}[2]{\textstyle{\frac {#1}{#2}}}
\newcommand{\ind}{\mathds}

\def\1{{\mathchoice {1\mskip-4mu\mathrm l}      
{1\mskip-4mu\mathrm l} 
{1\mskip-4.5mu\mathrm l} {1\mskip-5mu\mathrm l}}} 
 
\def\comment#1{} 
\newtheoremstyle{thm}{2ex}{2ex}{\itshape\rmfamily}{} 
{\bfseries\rmfamily}{}{1.7ex}{} 
 
\newtheoremstyle{rem}{1.3ex}{1.3ex}{\rmfamily}{} 
{\itshape\rmfamily}{}{1.5ex}{}

\renewcommand{\theequation}{\thesection.\arabic{equation}}

 \DeclareFontFamily{U}{mathx}{\hyphenchar\font45}
\DeclareFontShape{U}{mathx}{m}{n}{
      <5> <6> <7> <8> <9> <10>
      <10.95> <12> <14.4> <17.28> <20.74> <24.88>
      mathx10
      }{}
\DeclareSymbolFont{mathx}{U}{mathx}{m}{n}
\DeclareFontSubstitution{U}{mathx}{m}{n}
\DeclareMathSymbol{\bigtimes}{1}{mathx}{"91}

\renewcommand{\d}{{\rm d}} 
 
\newcommand{\eps}{\varepsilon}

\newcommand{\supp}{{\operatorname {supp}}} 
 
\newcommand{\dist}{{\operatorname {dist}}} 
\newcommand{\diam}{{\operatorname {diam}}}

\def\supp{\mathrm{supp}}

\newcommand{\Gcal}   {{\mathcal G }}

\newcommand{\Ocal}   {{\mathcal O }}

\newcommand\numberthis{\addtocounter{equation}{1}\tag{\theequation}}
\renewcommand{\e}   {{\operatorname e }}

\renewcommand{\cite}{\citet}

\theoremstyle{plain}
\newtheorem{theorem}{Theorem}[section]                                          
\newtheorem{proposition}[theorem]{Proposition}                          
\newtheorem{lemma}[theorem]{Lemma}
\newtheorem{corollary}[theorem]{Corollary}
\newtheorem{conjecture}[theorem]{Conjecture}
\theoremstyle{definition}
\newtheorem{definition}[theorem]{Definition}
\theoremstyle{remark}

\newtheorem{step}{STEP}

\makeatletter \@addtoreset{equation}{section} \makeatother
\renewcommand\theequation{\thesection.\arabic{equation}}
\renewcommand\thefigure{\thesection.\arabic{figure}}
\renewcommand\thetable{\thesection.\arabic{table}}

\newcommand{\aleaIndex}[1]{\href{http://alea.impa.br/english/index_v#1.htm}{\bf #1}}
\eheader{Alea}{\aleaIndex{7}}{2010}{21}{24}
\newcommand{\ProofEnde}{\hfill {$\square$}}

\renewcommand{\Z}{\mathbb{Z}}

\renewcommand{\R}{\mathbb{R}}

\def\SINR{\mathrm{SINR	}}

\begin{document}

\title[SINR percolation for Cox point processes]{Signal-to-interference ratio percolation \\ for Cox point processes}

\author{András Tóbiás}

\address{TU Berlin\newline
Straße des 17. Juni 136,\newline
10623 Berlin, Germany.}

\email{ tobias@math.tu-berlin.de}
\urladdr{\url{http://page.math.tu-berlin.de/~tobias/}}

\thanks{Research of the author supported by the Berlin Mathematical School}


\begin{abstract}
We study the signal-to-interference ratio (SINR) percolation model for a stationary Cox point process in two or higher dimensions, in case of a bounded and integrable path-loss function. We show that if this function has compact support or if the stationary intensity measure evaluated at a unit box has some exponential moments, then the SINR graph has an infinite connected component in case the spatial density of points is large enough and the interferences are sufficiently reduced (without vanishing). This holds under suitable stabilization and connectivity assumptions on the intensity measure. We also provide estimates on the critical interference cancellation factor.
\end{abstract}

\maketitle

\setcounter{tocdepth}{3}


\setcounter{section}{0}

\section{Introduction} \label{sec-Intro}
\subsection{Background and motivation}
Continuum percolation was introduced by Gilbert \cite{G61}. In his random graph model, two points of a homogeneous Poisson point process $X^\lambda$ in $\R^2$ with intensity $\lambda>0$ are connected by an edge if their distance is less than a fixed connection radius $r>0$. He showed that this model undergoes a phase transition: there is a critical intensity $\lambda_{\mathrm c}(r) \in (0,\infty)$ such that almost surely, for $\lambda<\lambda_{\mathrm c}(r)$ the graph consists of finite connected components, while for $\lambda>\lambda_{\mathrm c}(r)$ it \emph{percolates}, i.e., it has an infinite connected component. The motivation of this setting was to model a telecommunication network, in which the points of $X^\lambda$ are the users, and transmissions between users are only possible along the edges of the graph. In this view, long-distance communication is only possible if the graph percolates. 

The model of \cite{G61} has been widely studied and generalized in the Poisson case, see e.g.~\cite{MR96,FM07,BB09} for overviews. After 2010, it has also been extended to various other kinds of point processes, for example sub-Poisson \cite{BY10, BY13}, Ginibre and Gaussian zero \cite{GKP16}, and Gibbsian \cite{J16, S13}. The case of Gibbsian point processes was also studied earlier, see the references in \cite{J16}.

\cite{HJC17} considered Gilbert's graph model for a Cox point process, that is, a Poisson point process in a random environment. More precisely, let $\lambda>0$ and a stationary random measure $\Lambda$ on $\R^d$, $d \geq 2$, be given. A Cox point process $X^\lambda$ with intensity $\lambda \Lambda$ is characterized by the property that conditional on $\Lambda$, $X^\lambda$ is a Poisson point process with intensity $\lambda \Lambda$. In \cite{HJC17}, it was shown that under certain stabilization and connectedness conditions on $\Lambda$, $0<\lambda_{\mathrm c}(r)<\infty$ holds. More precisely, $\lambda_{\mathrm c}(r)>0$ if $\Lambda$ is stabilizing and $\lambda_{\mathrm c}(r)<\infty$ under the stronger assumption that $\Lambda$ is asymptotically essentially connected. These assertions are to be understood in the annealed sense, i.e., under a probability measure that governs $\Lambda$ and $X^\lambda$ jointly. 

According to \cite{HJC17}, the most important examples of $\Lambda$ for telecommunication are given by a stationary tessellation process, e.g., a Poisson--Voronoi, Poisson--Delaunay or Poisson line tessellation. The edge set of such a tessellation process can be used for modelling a telecommunication network on a street system, where the points of the Cox point process are the users, situated on the streets. The randomness of the tessellation process can be interpreted as the statistical variability of street systems in different areas. While Poisson--Delaunay tessellations fit well for modelling rural areas, Poisson--Voronoi tessellations are good approximations for various kinds of urban environments \cite{CGHJNP18}.

Another variant of Gilbert's graph model motivated by telecommunication is the signal-to-interference-plus-noise ratio (SINR) graph, which was considered in \cite{DBT03, DF06, FM07} in the case of a homogeneous Poisson point process with intensity $\lambda>0$ in $\R^2$. Here, two points are connected if the SINR between them is larger than a given threshold $\tau>0$ in both directions. The SINR of a transmission from $x \in \R^d$ to $y \in \R^d$ has the form $P\ell(|x-y|)/(N_0+\gamma P I(x,y))$. Here $\ell$ is a path-loss function describing the propagation of signal strength over distance, assumed decreasing, $P>0$ is the transmitted power, $N_0 \geq 0$ is the external noise and $I(x,y)$, the interference for the transmission from $x$ to $y$, is the sum of $\ell(|X_i-y|)$ over all Poisson points $X_i \notin \lbrace x, y \rbrace$, and $\gamma \geq 0$ is a factor expressing how strongly interference is cancelled compared to the signal. The motivation for the SINR model is that in real telecommunication networks, even if the transmitter is close to the receiver, the transmission may be unsuccessful due to too many other transmitters standing near the receiver, see~\cite[Section 1.2.5]{FM07}. Further, it is explained in \cite[Section~I]{DBT03} that $\gamma$ is a certain orthogonality factor between different transmissions, introduced by analogy with CDMA (i.e., Code Division Multiple Access) networks, which in applications usually takes values in $(0,1]$, and which appears due to the imperfect orthogonality of the codes used in CDMA. In order to simplify the presentation, throughout the rest of the paper we will assume that $P=1$.

For $\lambda>0$, let us write $\gamma^*(\lambda)$ for the supremum of all $\gamma>0$ for which the SINR graph percolates. If $\gamma=0$, then the SINR graph equals Gilbert's graph with radius $r_{\mathrm B}=\ell^{-1}(\tau N_0)$, and this graph contains all SINR graphs with positive $\gamma$. Thus, for $\lambda<\lambda_{\mathrm c}(r_{\mathrm B})$, we have $\gamma^*(\lambda)=0$. \cite{DF06} showed that under suitable integrability and boundedness assumptions on $\ell$, for any $\lambda>\lambda_{\mathrm c}(r_{\mathrm B})$, one has $\gamma^*(\lambda)>0$. Further, the following assertions were derived in \cite{DBT03, FM07} about $\lambda \mapsto \gamma^*(\lambda)$. SINR graphs with $\gamma>0$ have degrees bounded by $1+1/(\tau\gamma)$, which yields that $\gamma^*(\lambda) \leq 1/\tau$ for all $\lambda>0$. Further, $\gamma^*(\lambda)=\Ocal(1/\lambda)$ holds as $\lambda \to \infty$, and also $\gamma^*(\lambda)=\Omega(1/\lambda)$ if $\ell$ has bounded support. Here, for $f,g \colon (0,\infty) \to (0,\infty)$, we wrote $f(x)=\Ocal(g(x))$ if there exists $M,C>0$ such that $f(x) \leq Cg(x)$ for all $x>M$, and $f(x)=\Omega(g(x))$ if $g(x)=\Ocal(f(x))$. In \cite{BY13}, a more general notion of SINR graphs was considered, and the results of \cite{DF06} were extended to the case of sub-Poisson point processes in this context.

In particular, the upper bound of $\Ocal(1/\lambda)$ on the critical interference cancellation factor $\gamma^*(\lambda)$ for $\lambda$ large implies that adding more users to the telecommunication network (i.e., increasing $\lambda$) can actually harm percolation. This is a striking feature of the SINR graph model, while in e.g.~Gilbert's graph model increasing the intensity measure always improves the connectivity of the graph.

\subsection{Our findings} In the present paper, we investigate SINR percolation for Cox point processes, combining the benefits of modelling both user locations and connections between the users more realistically than in Gilbert's original model. To the best of our knowledge, this paper is the first one that considers SINR percolation also in $d \geq 3$ dimensions, despite the fact that some results of previous work about $d = 2$ extend to $d \geq 3$ without additional effort. We formulate our results for $d \geq 2$ whenever possible, and we point out which assertions of prior work extend to the higher dimensions. 

Let us summarize our most important results. First, we give general sufficient criteria for the existence of an infinite connected component in this model, in the case of a bounded path-loss function $\ell$. We consider the above defined SINR graph on a Cox point process $X^\lambda$ with intensity measure $\lambda\Lambda$ on $\R^d$, $d \geq 2$. We show that if $\Lambda$ is asymptotically essentially connected, then for $\lambda$ sufficiently large, we have $\gamma^*(\lambda)>0$ if any of the following additional assumptions is satisfied: (a) $\ell$ has bounded support, (b) for any compact set $A \subset \R^d$, there exists $\alpha>0$ such that $\mathbb E[\exp(\alpha \Lambda(A))]<\infty$, $\Lambda$ is $b$-dependent for some $b>0$, and $\ell$ satisfies the integrability assumption known from the Poisson case \cite{D71, DF06}, that is, $\int_{\R^d} \ell(|x|)\d x <\infty$. 

For the particular case of the homogeneous Poisson point process, our results imply that $\gamma^*(\lambda)>0$ holds for $\lambda$ sufficiently large, also for $d \geq 3$. The same holds for the Poisson--Voronoi and Poisson--Delaunay tessellations in any dimension in case $\ell$ has bounded support. Relevant examples of $\Lambda$ satisfying the condition (b), such that $\Lambda(A)$ is unbounded for bounded sets $A \subset \R^d$ with positive Lebesgue measure, include some intensity measures given by a shot-noise field with a compactly supported kernel function. 


Further, in the case when $d \geq 2$ and $\Lambda$ is only stabilizing, we show that if the connection radius $r_{\mathrm B}$ is sufficiently large, then $\gamma^*(\lambda)>0$ holds for sufficiently large $\lambda$ in case (b) above. This covers a number of relevant examples that are $b$-dependent but not asymptotically essentially connected.

We also provide estimates on $\gamma^*(\lambda)$. First, we conclude that the degree bound $1+1/(\tau\gamma)$ and the estimate that $\gamma^*(\lambda) \leq \frac{1}{\tau}$ for all $\lambda>0$ also hold in the Cox case. Second, we 
observe that if the number of Cox points who can successfully submit to a given point is bounded by some $k\in\N$ for all points, then every point can only receive messages from its $k$ nearest neighbours. This together with a high-confidence result of \cite{BB08} leads us to the conjecture that in the two-dimensional Poisson case, $\gamma^*(\lambda) \leq \frac{1}{4\tau}$ holds. Third, we show that for $b$-dependent Cox point processes, $\lim_{\lambda \to \infty} \gamma^*(\lambda)=0$. Here, our proof is applicable unless $N_0=0$ and $\ell$ has unbounded support. Finally, for $d=2$, for $b$-dependent Cox point processes with intensity measures that are locally bounded away from 0, we show that $\gamma^*(\lambda)=\Ocal(1/\lambda)$ holds as $\lambda \to \infty$, and also $\gamma^*(\lambda)=\Omega(1/\lambda)$ if also the support of $\ell$ is compact; these assertions generalize \cite[Theorem 4]{DBT03}.


The rest of this paper is organized as follows. In Section~\ref{sec-modeldefmainresults}, we define the model and present our main results. 
In particular, in Sections~\ref{sec-Gilbertmodeldef}, we summarize the results of \cite{HJC17} about continuum percolation for Cox point processes, and in Section~\ref{sec-SINRgraph} the ones of \cite{DBT03, DF06} about SINR percolation in the Poisson case. In Section~\ref{sec-phasetransitiondescription}, we present our main results about phase transitions in the Cox--SINR setting. Section~\ref{sec-gammalambda} contains our assertions and conjectures about the critical interference cancellation factor. In Section~\ref{sec-applicability} we discuss the applicability of the assertions of Sections~\ref{sec-phasetransitiondescription} and \ref{sec-gammalambda} to the main examples of the intensity measure. In Section~\ref{sec-phasetransitionproofs} we carry out and discuss the proofs of the results of Section~\ref{sec-phasetransitiondescription}. Finally, in Section~\ref{sec-gammalambdaproof} we verify the assertions of Section~\ref{sec-gammalambda}. 
 
\section{Model definition and main results}\label{sec-modeldefmainresults}
\subsection{Continuum percolation for Cox point processes}\label{sec-Gilbertmodeldef}\index{percolation!continuum!for Cox point processes}
In this section we recall the continuum percolation model defined in \cite[Section~2]{HJC17}. Let $\Lambda$ be a random element in the space $\mathbb M$ of\index{measure!Borel} Borel measures on $\R^d$, equipped with the\index{evaluation $\sigma$-field} evaluation $\sigma$-field \cite[Section 13.1]{LP17}, that is, the smallest $\sigma$-field that makes the mappings $B \mapsto \Lambda(B)$ measurable for all Borel sets $B \subseteq \R^d$. We always assume that $d \geq 2$. 
We define 
\[ Q_r(x)=x+[-r/2,r/2]^d \]
for $x \in \R^d$ and $r \geq 0$, further, we write $Q_r=Q_r(o)$, where $o$ denotes the origin of $\R^d$.  We assume that $\mathbb E[\Lambda(Q_1)]=1$ and $\Lambda$ is\index{measure!intensity!stationary}\index{point process!stationary} stationary, that is, $\Lambda(\cdot)$ equals $\Lambda(\cdot + x)$ in distribution for all $x \in \R^d$. 

Then for $\lambda>0$, we let $X^\lambda$ be a\index{point process!Cox|textbf} Cox point process with intensity $\lambda \Lambda$. That is, conditional on $\Lambda$, $X^\lambda$ is a Poisson point process with intensity\index{measure!intensity} $\lambda \Lambda$.
Note that the conditions on $\Lambda$ imply that for all $x \in \R^d$, $\Lambda(\lbrace x\rbrace)=0$ holds almost surely, and thus $X^\lambda$ is a\index{point process!simple} \emph{simple} point process. That is, one can write $X^\lambda=(X_i)_{i \in \N}$, so that almost surely, for all $i,j \in \N$, $X_i \neq X_j$ unless $i=j$. Further, if $\Lambda$ is identically equal to the Lebesgue measure $ |\cdot|$, then $X^\lambda$ is a homogeneous Poisson point process with intensity $\lambda$. We will often simply say ``Cox process'' instead of ``Cox point process''. We denote by $\Lambda_B$ the restriction of the random measure $\Lambda$ to the set $B \subset \R^d$. 

Let us give some examples of random intensity measures satisfying our assumptions. Any\index{measure!intensity!absolutely continuous} absolutely continuous example $\Lambda$ has the form $\Lambda(\d x)=l_x \d x$ for a stationary non-negative\index{random field!stationary} random field $l=\lbrace l_x \rbrace_{x \in \R^d}$ with $\mathbb E[l_o]=1$, see \cite[Example 2.1]{HJC17}. Examples include the\index{point process!Poisson, modulated} modulated Poisson point process: $l_x = \lambda_1 \ind{1}\lbrace x \in \Xi \rbrace+\lambda_2 \ind{1}\lbrace x \notin \Xi \rbrace$ for a stationary\index{random closed set} random closed set $\Xi$  and $\lambda_1,\lambda_2 \geq 0$, and intensities given by a\index{shot-noise field} shot-noise field: $l_x=\sum_{Y_i \in Y_{\mathbf S}} k(x-Y_i)$ for a non-negative integrable kernel $k \colon \R^d \to [0,\infty)$ with compact support and $Y_{\mathbf S}$ a Poisson point process with intensity $\lambda_{\mathbf S}>0$. Relevant singular\index{measure!intensity!singular} examples are the Poisson point processes on\index{street system!random} random street systems \cite[Example 2.2]{HJC17}. Here, $\Lambda(\d x) = \nu_1(S \cap \d x)$ for a stationary point process $S$ with values in the space of line segments, e.g., a\index{tessellation!Poisson--Voronoi} Poisson--Voronoi or \index{tessellation!Poisson--Delaunay} Poisson--Delaunay tessellation, where $\nu_1$ denotes\index{measure!Hausdorff!one-dimensional} one-dimensional Hausdorff measure.

For $r,\lambda>0$, the\index{graph!Gilbert|textbf}\index{graph!Gilbert!Cox} \emph{Gilbert graph} $g_r(X^\lambda)$ is defined as follows. Its vertex set is $X^\lambda$, more precisely the set $\lbrace X_i \colon i \in \N \rbrace$, and $X_i,X_j \in X^\lambda$, $i \neq j$,  are connected by an edge whenever their distance is less than the connection radius $r$. A\index{cluster} \emph{cluster} in a (possibly random) graph is a maximal connected component, and we say that\index{percolation!in the Gilbert graph} the graph \emph{percolates} if it contains an\index{cluster!infinite} infinite cluster. The\index{critical intensity!for the Gilbert graph}\index{phase transition} critical intensity is defined as
\[ \lambda_{\mathrm c}(r) = \inf \lbrace  \lambda  >0 \colon \mathbb P(g_r(X^\lambda)\text{ percolates})>0 \rbrace. \]
Percolation of $g_r(X^\lambda)$ occurs if and only if the\index{percolation!in the Boolean model} associated\index{Boolean model|textbf} \emph{Boolean model}, that is, $X^\lambda \oplus B_{r/2} = \bigcup_{i \in \N} B_{r/2}(X_i)$, has an unbounded connected component, see \cite[Section~7.1]{HJC17}. Here we wrote $B_R=B_R(o)$ where $B_R(x)$ denotes the open $\ell^2$-ball of radius $R$ around $x$ for $x \in \R^d$ and $R>0$. Note that for fixed $r>0$, $\lambda \mapsto \mathbb P(g_r(X^\lambda)\text{ percolates})$ is increasing in $\lambda$. Given $r>0$, any intensity $\lambda \in (0,\lambda_{\mathrm c}(r))$ is called\index{phase!subcritical} \emph{subcritical}, $\lambda=\lambda_{\mathrm c}(r)$\index{phase!critical} \emph{critical}, and any $\lambda \in (\lambda_{\mathrm c}(r),\infty)$\index{phase!supercritical} \emph{supercritical}. 

The next two definitions are crucial in \cite{HJC17} for showing that a subcritical respectively supercritical phase exists. The first notion is\index{stabilization|textbf} stabilization, which means a certain decay of spatial correlations of the intensity measure with distance.  We let $\dist_p(\varphi,\psi)=\inf \lbrace \Vert x-y \Vert_p \colon x\in\varphi, y \in \psi \rbrace$ denote the $\ell^p$-distance between two sets $\varphi,\psi \subset \R^d$ for $p \in [1,\infty]$.
\begin{definition}\label{defn-stabilization}
The random measure $\Lambda$ is \emph{stabilizing} if there exists a random field of\index{stabilization!radii} \emph{stabilization radii} $R=\lbrace R_x \rbrace_{x\in \R^d}$ defined on the same probability space as $\Lambda$ such that, writing
\[ R(Q_n(x)) = \sup_{y \in Q_n(x) \cap \Q^d} R_{y} ,\qquad n \geq 1, ~x \in \R^d, \]
the following hold.
\begin{enumerate}
\item\label{first-stabilizing} $R$ is measurable with respect to $\Lambda$, and $(\Lambda,R)$ is jointly stationary,
\item $\lim_{n \to \infty} \mathbb P(R(Q_n)<n)=1$,
\item for all $n \geq 1$, for any bounded measurable function $f \colon \mathbb M \to [0,\infty)$ and finite $\varphi \subseteq \R^d$ with $\dist_2(x,\varphi\setminus \lbrace x \rbrace)>3n$ for all $x \in \varphi$, the following random variables are independent:
\[  f(\Lambda_{Q_n(x)}) \ind{1}\lbrace R(Q_n(x)) <n \rbrace, \qquad x \in \varphi. \] 
\end{enumerate}
\end{definition}

A strong form of stabilization is \emph{$b$-dependence}; for $b>0$, $\Lambda$ is called $b$-dependent\index{b-dependence@$b$-dependence|textbf} if $\Lambda_A$ and $\Lambda_B$ are independent whenever $\dist_2(A,B)>b$. On the other hand, in this paper, $b$-dependence of stochastic processes defined on discrete subsets of $\R^d$ with an explicitly given value of $b$ will always be meant with $\dist_\infty$ instead of $\dist_2$ on the discrete set; we note that \cite{HJC17} also used this convention tacitly.

Let us write
$\supp(\mu) = \lbrace x \in \R^d \colon \mu(Q_{\eps}(x))>0,~\forall \eps>0 \rbrace$
for the\index{measure!intensity!support}\index{support!of a (possibly singular) measure} support of a (possibly singular) measure $\mu$. 
The second notion, asymptotic essential connectedness, indicates, in addition to stabilization, strong local connectivity of the intensity measure. 
\begin{definition}\label{defn-asessconn}
The stabilizing random measure $\Lambda$ with stabilization radii $R$ is\index{asymptotic essential connectedness|textbf} \emph{asymptotically essentially connected} if for all sufficiently large $n \geq 1$, whenever $R(Q_{2n})<n/2$, we have that
\begin{enumerate}
    \item $\supp(\Lambda_{Q_n})$ contains a connected component of diameter at least $n/3$, and
    \item if $C$ and $C'$ are connected components in $\supp(\Lambda_{Q_n})$ of diameter at least $n/9$, then they are both contained in one of the connected components of $\supp(\Lambda_{Q_{2n}})$.
\end{enumerate}
\end{definition}

As for the main examples,\index{stabilization!examples} it was shown in \cite[Section~3.1]{HJC17} that Poisson--Voronoi and Poisson--Delaunay tessellations are\index{asymptotic essential connectedness!examples} asymptotically essentially connected. Further, shot-noise fields are\index{b-dependence@$b$-dependence!examples} $b$-dependent but only in some cases asymptotically essentially connected, see Section~\ref{sec-examplesphasetransition} for further details. The modulated Poisson point process is also $b$-dependent if $\Xi$ is a\index{Boolean model!Poisson}\index{point process!Poisson, modulated!by a Poisson--Boolean model} Poisson--Boolean model (that is, the Boolean model of a homogeneous Poisson point process). In this case, it is also asymptotically essentially connected if $\lambda_1,\lambda_2>0$, or if $\lambda_1>\lambda_2=0$ and the underlying Poisson--Boolean model $\Xi$ is supercritical. However, in general both for $\lambda_1>\lambda_2=0$ and for $\lambda_2>\lambda_1=0$ it may happen that $\Lambda$ is not asymptotically essentially connected, as we will explain in Section~\ref{sec-examplesphasetransition}.  
Poisson line tessellations and their rectangular variants like Manhattan grids are also relevant for modelling street systems \cite{GFSS05, HHJC19}, however, they are not stabilizing, and neither the existence of subcritical phase nor the one of supercritical phase has been verified for them so far.

By \cite[Theorems 2.4, 2.6]{HJC17}, for any $r>0$ the following holds. If $\Lambda$ is stabilizing, then $\lambda_{\mathrm c}(r)>0$. If $\Lambda$ is asymptotically essentially connected, then $\lambda_{\mathrm c}(r)<\infty$.

In these results, roughly speaking, the spatial decorrelation coming from stabilization of $\Lambda$ makes it easy to verify, using discrete percolation techniques, that long-distance connections in $g_r(X^\lambda)$ do not exist for $\lambda>0$ sufficiently small, see~\cite[Section 5.1]{HJC17}. On the other hand, as $\lambda \to \infty$, $X^\lambda$ fills the support of $\Lambda$ with high probability. This fact together with the stabilization of $\Lambda$ and the strong connectivity of the support of $\Lambda$ can be used in order to verify percolation of $g_r(X^\lambda)$ for large $\lambda$ if $\Lambda$ is asymptotically essentially connected, cf.~\cite[Section 5.2]{HJC17}.

\subsection{Signal to interference plus noise ratio graph}\label{sec-SINRgraph}\index{graph!SINR|textbf} In this section we follow \cite{DF06}. We choose a decreasing\index{path-loss function} \emph{path-loss function} $\ell: [0,\infty) \to [0,\infty)$, which describes the propagation of signal strength over distance. Note that $\ell(|x-y|) \leq \ell(0)$ for all $x,y \in \R^d$. Further assumptions on $\ell$ will be made below using the following definitions. For two points $X_i,X_j$ of the Cox point process $X^\lambda$, we define the \emph{signal-to-interference-plus noise ratio (SINR)} of the transmission from $X_i$ to $X_j$ as follows\index{SINR|textbf}
\[  \SINR(X_i,~X_j,~X^\lambda) = \frac{\ell(|X_i-X_j|)}{N_0 + \gamma \sum_{k \neq i,j} \ell(|X_k-X_j|)}, \numberthis\label{SINR}\]
where $N_0 \geq 0$ is the environmental\index{noise} \emph{noise}, the sum in the denominator of \eqref{SINR} is called the\index{interference|textbf} \emph{interference} (of the transmission from $X_i$ to $X_j$), and $\gamma \geq 0$ is the\index{interference cancellation factor} \emph{interference cancellation factor}.
Then we fix\index{SINR!threshold} $\tau>0$ and say that the transmission from $X_i$ to $X_j$ is possible if and only if
\[ \SINR(X_i,~X_j,~X^\lambda) >\tau. \numberthis\label{SINR>tau} \]

We will tacitly exclude the\index{graph!SINR!with $\gamma=N_0=0$ is degenerate} degenerate case $\gamma=N_0=0$. Further, if $N_0=0$, we use the\index{SINR!threshold!convention for $N_0=0$} convention \cite[Section~6.1]{BB09} that the inequality \eqref{SINR>tau} holds if $\ell(|X_i-X_j|)>\tau \gamma \sum_{k \neq i,j} \ell(|X_k-X_j|)$.

We define the\index{graph!SINR!directed} \emph{directed SINR graph} $g^{\rightarrow}_{(\gamma,N_0,\tau)}(X^\lambda)$ on the vertex set $X^\lambda$ via drawing a directed edge pointing from $X_i$ towards $X_j $ (denoted as $X_i \to X_j$) whenever $i \neq j$ and $\SINR(X_i,~X_j,~X^\lambda) >\tau$. Next, the\index{graph!SINR!undirected} (undirected) \emph{SINR graph} $g_{(\gamma,N_0,\tau)}(X^\lambda)$ has vertex set $X^\lambda$, and $(X_i,X_j) \in X^{\lambda} \times X^{\lambda}$ is an edge in $g_{(\gamma,N_0,\tau)}(X^\lambda)$ if and only if both $X_i \to X_j$ and $X_j \to X_i$ are edges in $g^{\rightarrow}_{(\gamma,N_0,\tau)}(X^\lambda)$. 

We note that \cite{KY07} studied percolation in the directed SINR graph in the two-dimensional Poisson case. It obtained results that are very similar to the ones of \cite{DBT03, DF06, FM07} in the undirected case. In the present paper we will focus on the undirected SINR graph, but we will also use some properties of the directed one in our arguments.

See Figure~\ref{figure-gammas} for simulations of the SINR graph in the two-dimensional Poisson case. Extensions of the SINR graph model, such as external interferers \cite{BY13}, the information theoretically secure SINR graph \cite{VI14}, and random signal powers, were surveyed in the author's PhD thesis \cite[Sections 4.2.3.4--4.2.3.6]{T19}, including a number of remarks about the case of Cox point processes.


As for $N_0>0$ and $\gamma=0$, $X_i,X_j$ are connected by an edge in $g_{(0,N_0,\tau)}(X^\lambda)$ if and only if the\index{SNR} \emph{signal-to-noise ratio (SNR)} between them is larger than $\tau$, i.e.,
\[ \mathrm{SNR}(X_i,~X_j)=\mathrm{SNR}(X_j,~X_i)=\frac{\ell(|X_i-X_j|)}{N_0}>\tau. \numberthis\label{SNRcondition} \]
Whenever $\ell^{-1}(\tau N_0)$ is well-defined and positive (in particular, $\ell(0)>\tau N_0$), this is equivalent to
$ |X_i-X_j| \leq \ell^{-1}( \tau N_0) . $
In this case\index{graph!SINR!for $\gamma=0$}\index{graph!SNR} $g_{(0,N_0,\tau)}(X^\lambda)$ equals the Gilbert graph\index{SNR!threshold} $g_{r_{\mathrm B}}(X^\lambda)$, where
\[ r_{\mathrm B} = \ell^{-1}(\tau N_0) . \numberthis\label{criticalradius} \]
For two graphs $G=(V,E),G'=(V,E')$ with the same vertex set $V$, we will write $G \preceq G'$ if $E \subseteq E'$, i.e., if all edges in $G$ are also contained in $G'$. Now, for $\tau>0$ and $N_0 > 0$, we have $g_{(\gamma,N_0,\tau)}(X^\lambda) \preceq g_{(\gamma',N_0,\tau)}(X^\lambda)$ for all $0 \leq \gamma'<\gamma$. Thus,\index{graph!SINR!decreasing in $\gamma$} $g_{(\gamma,N_0,\tau)}(X^\lambda) \preceq g_{r_{\mathrm B}}(X^\lambda)$, hence any edge of $g_{(\gamma,N_0,\tau)}(X^\lambda)$ has length at most $r_{\mathrm B}$. In contrast, if $N_0=0$ (and $\gamma>0$), then the edge lengths of the SINR graph $g_{(\gamma,0,\tau)}(X^\lambda)$ are unbounded. On the other hand, while Gilbert graphs have no bound on the degrees of the vertices, we will show in Section~\ref{sec-gammalambda} that\index{degree bounds} all in-degrees in $g^{\rightarrow}_{(\gamma,N_0,\tau)}(X^\lambda)$ and thus also all degrees in $g_{(\gamma,N_0,\tau)}(X^\lambda)$ are bounded by $1+1/(\tau\gamma)$ for fixed $\gamma>0$. If additionally also $N_0>0$, then an easy computation of SINR values implies
that points that are not isolated in $g_{(\gamma,N_0,\tau)}(X^\lambda)$ have uniformly bounded degrees even in the Gilbert graph $g_{r_{\rm B}}(X^\lambda)$. 

In Section~\ref{sec-generalgammabound} we will explain that for stationary Cox processes, it follows easily that whenever $\gamma>0$ and in-degrees in $g^{\rightarrow}_{(\gamma,N_0,\tau)}(X^\lambda)$ are bounded by $k \in \N$, any two points that are connected in  $g_{(\gamma,N_0,\tau)}(X^\lambda)$ are mutually among the $k$ nearest neighbours of each other in $X^\lambda$. This assertion was not explicitly mentioned in earlier works about SINR percolation, even though its proof is immediate. This assertion implies that $g_{(\gamma,N_0,\tau)}(X^\lambda)$ is a subgraph of the \emph{bidirectional $k$-nearest neighbour graph} considered in \cite{BB08}, where two points $X_i \neq X_j$ of $X^\lambda$ are connected by an edge whenever $X_i$ is one of the $k$ nearest neighbours of $X_j$ in $X^\lambda$ and also vice versa. We will elaborate on some possible consequences of this relation in Section~\ref{sec-gammalambda}.

Now, we define\index{interference cancellation factor!critical}\index{percolation!in the SINR graph}
\[  \gamma^*(\lambda)=\gamma^*(\lambda,N_0,\tau):=\sup \lbrace \gamma>0 \colon \P \big(g_{(\gamma,N_0,\tau)}(X^\lambda) \text{ percolates}\big)>0  \rbrace \numberthis\label{gammastarlambda} \]
for fixed $\lambda,\tau>0$ and $N_0 \geq 0$. Further, we put\index{critical intensity!for the SINR graph}
\[ \lambda_{N_0,\tau}= \inf \lbrace \lambda>0 \colon 
\gamma^*(\lambda')>0,~\forall \lambda' \geq \lambda 
\rbrace. \numberthis\label{SINRlambdac} \]
Then for $\lambda<\lambda_{N_0,\tau}$, $\mathbb P(g_{(\gamma,N_0,\tau)}(X^\lambda)\text{ percolates})=0$ for all $\gamma>0$. \emph{A priori}, there is no reason for $\lambda_{N_0,\tau}=\inf \lbrace \lambda>0 \colon \gamma^*(\lambda)>0 \rbrace$ to hold, but this identity will turn out to be true in most of the cases that we consider.

For any Cox--SINR graph with $\gamma \geq 0$, the existence of an infinite cluster is a shift-invariant event. Therefore,\index{cluster!infinite!probability of existence} if the stationary intensity measure $\Lambda$ is also ergodic, then the probability of this event is either zero or one, and the\index{cluster!infinite!number of infinite clusters} number of infinite clusters is almost surely constant (possibly infinite), cf.~\cite[Theorem~2.1]{MR96}. In particular, this holds for stabilizing Cox processes, since it is easy to derive that\index{stabilization!implies mixing}\index{stabilization!implies mixing} stabilization implies mixing and therefore also\index{stabilization!implies ergodicity} ergodicity. Without ergodicity, one can find examples where this property fails, see e.g.~\cite[Section 4.2.3.3]{T19}.

Now, let us fix $N_0\geq 0$, $\tau>0$, and let us make the following assumption on the path-loss function $\ell$ (which also implies that $\ell$ is decreasing) for the rest of this paper.

\noindent \textbf{Assumption ($\ell$)}.\index{Assumption ($\ell$)}
\begin{enumerate}[(i)]
\item\label{first-pathloss} $\ell$ is\index{path-loss function!continuity}\index{path-loss function!constant part at zero} continuous, constant on $[0,\upsilon_0]$ for some $\upsilon_0 \geq 0$, and on $[\upsilon_0,\infty) \cap \supp~ \ell$ it is strictly decreasing,
\item\label{third-pathloss} $1 \geq \ell (0) > \tau N_0$,
\item\label{last-pathloss}\index{path-loss function!integrability}
$\int_{\R^d} \ell(|x|) \d x <\infty$.
\end{enumerate}
These constraints on $\ell$ are slightly more general than the ones of \cite{DF06} because we allow $\upsilon_0$ to be positive, motivated by the facts that the proof of the main results of \cite{DF06} works also for $\upsilon_0>0$ and path-loss functions with $\upsilon_0>0$ are widely used in practice (see e.g.~the simulations in \cite{DBT03, DF06}). E.g., the path-loss function $\ell(r)=\min \lbrace 1, r^{-\alpha} \rbrace$, $\alpha>d$ (where we recall that $d$ is the dimension), corresponding to\index{Hertzian propagation} ideal Hertzian propagation \cite{DBT03} satisfies Assumption~($\ell$). 

Let us recall \cite[Theorem~1]{DF06} about the\index{point process!Poisson} homogeneous Poisson case $\Lambda \equiv |\cdot|$ for $d=2$.
\begin{theorem}[\cite{DF06}]\label{theorem-fromDF06}\index{percolation!SINR!for Poisson point processes}\index{percolation!SINR!in two dimensions}\index{phase transition!in the SINR graph!for Poisson point processes|textbf}
If $\Lambda \equiv |\cdot|$, $d=2$, and $N_0,\tau>0$, then $\lambda_{N_0,\tau}=\lambda_{\mathrm c}(r_{\mathrm B}) \in (0,\infty)$.
\end{theorem}
In words, for any intensity $\lambda$ such that the SNR graph $g_{(0,N_0,\tau)}(X^\lambda)=g_{r_{\mathrm B}}(X^\lambda)$ is supercritical, there exists a small but positive $\gamma$ such that $g_{(\gamma,N_0,\tau)}(X^\lambda)$ still percolates. The case $N_0=0$ will be discussed in Section~\ref{sec-noIse} in the general Cox case. Simulations of the SINR graph of a two-dimensional Poisson point process can be seen in Figure~\ref{figure-gammas}. 

 According to the results of \cite{D71}, for bounded path-loss functions $\ell$ not satisfying Assumption $(\ell)$ \eqref{last-pathloss}, the SINR graph of a Poisson point process has no edges for $\gamma>0$. However, this does not exclude percolation in the SINR graph in the case of an unbounded path-loss function $\ell \colon (0,\infty) \to [0,\infty)$ satisfying $\int_{\R^d \setminus B_\eps} \ell(|x|)\d x <\infty$ for all $\eps>0$. In fact, \cite{DBT03} conjectured that percolation occurs for $\ell(r)=r^{-\alpha}$, $\alpha>d$, in the two-dimensional Poisson case. Further, \cite[Section 3.3.2]{D05} has shown that in case of this path-loss function, if $\lambda_{N_0,\tau}>0$ holds true, then $\lambda \mapsto \gamma^*(\lambda)$ is increasing (but bounded thanks to the degree bounds, see Section~\ref{sec-generalgammabound}). In contrast, for path-loss functions satisfying Assumption $(\ell)$, $\gamma^*(\lambda)$ tends to zero as $\lambda \to \infty$, see Section~\ref{sec-lambdadependentgammabound}. See Figure~\ref{figure-gammastarlambda} for a visual sketch of the already verified and the conjectured properties of the function $\lambda \mapsto \gamma^*(\lambda)$ in the Cox case, which summarizes some of the main results of our paper.
\subsection{Phase transitions}\label{sec-phasetransitiondescription}\index{phase transition}
This section contains our main results about percolation properties of $g_{(\gamma,N_0,\tau)}(X^\lambda)$ depending on the parameters $N_0,\tau,\lambda,\gamma$. In Section~\ref{sec-positivenoise} we present our main results for fixed $N_0,\tau>0$. In this setting, the SNR radius $r_{\mathrm B}$ is fixed, and thus, according to \cite[Theorem 2.6]{HJC17}, if $\Lambda$ is asymptotically essentially connected, then the SNR graph percolates for large $\lambda$ with positive probability (actually with probability 1 thanks to stabilization). We show that under additional assumptions on $\Lambda$ and $\ell$, we have $\gamma^*(\lambda)>0$ for all sufficiently large $\lambda$. Under similar assumptions, in Section~\ref{sec-stabilizing} we show that if $\Lambda$ is only stabilizing, then one can choose the SNR radius $r_{\mathrm B}$ large enough  such  that $\gamma^*(\lambda)>0$ occurs for sufficiently large $\lambda>0$. In Section~\ref{sec-noIse} we comment on the case $N_0=0$. 

\subsubsection{The case of asymptotically essentially connected intensity}\label{sec-positivenoise}
If $\Lambda$ is\index{asymptotic essential connectedness} asymptotically essentially connected, then the SINR graph percolates for large enough $\lambda$ and accordingly chosen small enough $\gamma>0$ under additional assumptions on $\ell$ and $\Lambda$.\index{percolation!SINR!for Cox point processes|textbf}
\begin{theorem}\label{prop-firstpercolation} Let $N_0,\tau>0$.
\begin{enumerate}[(1)]\index{phase transition!in the SINR graph!for Cox point processes}
\item\label{firstsubcritical} $\lambda_{N_0,\tau} \geq \lambda_{\mathrm c}(r_{\mathrm B})$. Further, if $\Lambda$ is stabilizing, then $\lambda_{N_0,\tau} >0$.
\item\label{firstsupercritical} If $\Lambda$ is asymptotically essentially connected, then $\lambda_{N_0,\tau}<\infty$ holds if at least one of the following two conditions is satisfied:
\begin{enumerate}[(a)]
\item\label{ellboundedcond}\index{path-loss function!compactly supported}\index{support!of the path-loss function} $\ell$ has compact support,
\item\label{expmomentscond}\index{exponential moments of the intensity}\index{path-loss function!stronger decay condition} $\Lambda$ is $b$-dependent, and $\mathbb E[\exp(\alpha \Lambda(Q_1))]<\infty$ for some $\alpha>0$.
\end{enumerate}
\end{enumerate}
\end{theorem}
We already see that Theorem~\ref{prop-firstpercolation}\eqref{firstsubcritical} is true. Indeed, it follows from \cite[Theorem 2.4]{HJC17} and the fact that for $N_0,\tau,\gamma>0$, we have $g_{(\gamma,N_0,\tau)}(X^\lambda) \preceq g_{r_\mathrm B}(X^\lambda)$. Note that this assertion requires only that SNR graph be a well-defined Gilbert graph; for this, it suffices if $\ell: (0,\infty) \to [0,\infty)$ is decreasing and the radius $r_{\rm B}=\ell^{-1}(\tau N_0)$ is well-defined and positive. In particular, $\lim_{r \downarrow 0} \ell(r)=\infty$ is not a problem for this assertion. On the other hand, as we saw in Section~\ref{sec-SINRgraph}, unless $\ell$ has integrable tails, SINR graphs with $\gamma>0$ have no edges in the Poisson case. It is easy to see that the same holds in the general stationary Cox case. Thus, the integrability condition \eqref{last-pathloss} of Assumption $(\ell)$ in Theorem~\ref{prop-firstpercolation}\eqref{firstsupercritical} is optimal for percolation in the Cox--SINR graph in case of a bounded path-loss function. 

The proof of Theorem~\ref{prop-firstpercolation}\eqref{firstsupercritical} is carried out in Section~\ref{sec-firstpercolationproof}. The proof consists of four steps. First, we map our percolation problem to a discrete site percolation model. Second, we argue why this discrete model has an unbounded cluster for large $\lambda$ and accordingly chosen small $\gamma>0$, conditional on the assumption that interferences can be suitably controlled. Third, we show that if the discrete model percolates, then so does $g_{(\gamma,N_0,\tau)}(X^\lambda)$. Fourth, we finish the proof of percolation in the discrete model by controlling the interferences. At this point we use the assumption \eqref{ellboundedcond} respectively \eqref{expmomentscond}. In Section~\ref{sec-additionalchallenges} we will comment on possible generalizations of the proof, using the notation introduced in Section~\ref{sec-firstpercolationproof}.

In the rest of this paper, we will write ``$Y$ is bounded (away from 0)'' equivalently to ``$Y$ is almost surely bounded (away from 0)'' for any nonnegative random variable $Y$. Then we see that any $b$-dependent $\Lambda$ such that $\Lambda(Q_1)$ is bounded satisfies \eqref{expmomentscond}; an asymptotically essentially connected example is the modulated Poisson point process with $\Xi$ being a Poisson--Boolean model in case $\lambda_1,\lambda_2>0$, or also in case $\lambda_1>\lambda_2=0$ if $\Xi$ is supercritical.
As for unbounded intensity measures $\Lambda$ satisfying the exponential moment condition in \eqref{expmomentscond}, our main $b$-dependent examples are the shot-noise fields that are asymptotically essentially connected (see Section~\ref{sec-examplesphasetransition} for explicit examples). For these intensity measures, Theorem~\ref{prop-firstpercolation} implies that for large enough $\lambda$ and small enough $\gamma$ its SINR graph percolates in the case of any path-loss function $\ell$ satisfying Assumption $(\ell)$, in particular, for all $\ell$ satisfying \eqref{first-pathloss} and \eqref{third-pathloss} of that assumption and $\ell(r)=\Ocal(1/r^{2+\eps})$ for some $\eps>0$. 

B.~Jahnel and the author recently showed \cite{JT19} that the total edge length of two-dimensional Poisson--Voronoi or Poisson--Delaunay tessellations in a unit square has all exponential moments. As already mentioned, these tessellations are also asymptotically essentially connected, however, they are not $b$-dependent. Hence, only the condition \eqref{ellboundedcond} is applicable for them. Since these tessellations
can well be used for modelling statistical properties of real street systems \cite{GFSS05, CGHJNP18}, it is a highly interesting open question to verify percolation in the SINR graph for the Cox process on these tessellations in case of an unboundedly supported path-loss function. In higher dimensions, even the existence of exponential moments of $\Lambda(Q_1)$ is open. See Figure~\ref{figure-PVTgammas} for simulations of SINR graphs in the case of the Poisson--Voronoi tessellation for $d=2$. 



We see that Theorem~\ref{prop-firstpercolation}\eqref{firstsupercritical} holds\index{percolation!SINR!in higher dimensions} in particular in the\index{point process!Poisson} Poisson case and thus it generalizes Theorem~\ref{theorem-fromDF06} to $d \geq 3$ dimensions. However, it does not recover the identity that $\lambda_{N_0,\tau}=\lambda_{\mathrm c}(r_{\mathrm B})$ for all $N_0,\tau>0$. 
We nevertheless expect that the statement is also true in the higher-dimensional Poisson case. We defer the proof of this conjecture to future work. Further, in Section~\ref{sec-examplesphasetransition} we will discuss the applicability of all results of Section~\ref{sec-phasetransitiondescription} to each of the examples introduced in Section~\ref{sec-Gilbertmodeldef}.

In fact, $\lambda_{N_0,\tau}<\infty$ can also be proven for certain classes of non-stabilizing Cox processes, cf.~\cite[Sections 4.2.3.2, 4.2.3.3]{T19}. In some cases, even $\lambda_{N_0,\tau}=0$ holds \cite[Section 4.2.3.3]{T19}, which is impossible in the stabilizing case thanks to \cite[Theorem 2.4]{HJC17}.

Apart from the value of positive results on percolation in an SINR graph with $\gamma>0$, such as Theorem~\ref{prop-firstpercolation}\eqref{firstsupercritical}, for applications in telecommunications, such assertions have important theoretical consequences for more well-known continuum percolation models such as Gilbert and $k$-nearest neighbour graphs. Namely, the underlying Gilbert graph keeps percolating after removing all vertices that have degree larger than $n$, given that $n$ is large enough. On the other hand, the bidirectional $k$-nearest neighbour graph containing the SINR graph keeps percolating after removing all edges that have length larger than $r_{\mathrm B}$. This also implies the same statement for the more frequently studied \emph{undirected} $k$-nearest neighbour graph \cite{HM96,BB08}, where one connects two points whenever at least one of them is one of the $k$ nearest neighbours of the other. 

\subsubsection{The case of only stabilizing intensity}\label{sec-stabilizing}
\begin{figure}
    \centering
    \includegraphics[scale=1.9]{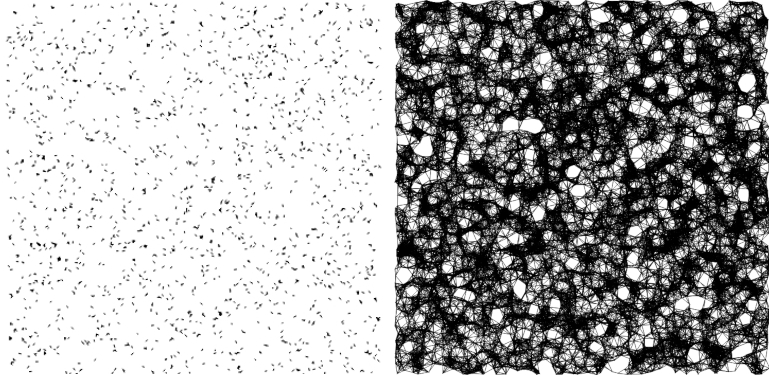}
    \caption{Gilbert graphs in case of a stabilizing intensity measure such that $\lambda_{\rm c}(r)<\infty$ holds only for large $r$. The intensity measure is given as $\lambda \Lambda(\d x)=\lambda \lambda_1 \ind{1}\lbrace x \in \Xi \rbrace \d x$, where $\lambda>0$ is very large, $\lambda_1>0$, and $\Xi$ is a strongly subcritical Poisson--Boolean model. (1) Small $r$: even though the density of Cox points per unit volume is high, the intensity measure has large void spaces. The Gilbert graph is split into many small connected components thanks to the disconnectedness of the support of $\Lambda$. (2) Large $r$: the small components can now connect up so that $g_r(X^\lambda)$ percolates.}
    \label{figure-MPP}
\end{figure}

\index{stabilization}
According to \cite[Section~2.1]{HJC17},\index{phase!supercritical!examples for absence for stabilizing $\Lambda$ and small $r$} stabilization of $\Lambda$ does not imply that $\lambda_{\mathrm c}(r)<\infty$ for all $r>0$, see Figure~\ref{figure-MPP} and Section~\ref{sec-examplesphasetransition} for more details. Now we argue that if $\Lambda$ is stabilizing with $\E[\Lambda(Q_1)]=1$, then $\lambda_{\mathrm c}(r)<\infty$ holds for $r$ large enough, and for the SINR graph, if $r_{\mathrm B}$ is large (in particular $N_0>0$) and the condition \eqref{expmomentscond} holds, then also $\lambda_{N_0,\tau}<\infty$. 

The fact that $\lambda_{\mathrm c}(r)<\infty$ holds for $r$ large for $\Lambda$ stabilizing is actually a direct consequence of certain results of \cite{HJC17}, but since it was not stated explicitly in that paper, we present it as a corollary.\index{stabilization!implies existence of supercritical phase for large $r$}
\begin{corollary}\label{cor-largersupercritical}\index{phase!supercritical!existence for stabilizing $\Lambda$ for large $r$}
If $\Lambda$ is stabilizing, then the following hold.
\begin{enumerate}
\item\label{first-largersupercritical}
There exists $r_0 \geq 0$ such that $\lambda_{\mathrm c}(r)<\infty$ holds for all $r > r_0$. \item\label{second-largersupercritical} $\lim_{r \to \infty} \lambda_{\mathrm c}(r)=0$.
\end{enumerate}
\end{corollary}
The proof of Corollary~\ref{cor-largersupercritical} is carried out in Section~\ref{sec-largersupercritical}. We will see that after recalling some elements of Palm calculus and the notion of percolation probability for Cox processes from \cite{HJC17}, the corollary follows immediately from \cite[Theorem~2.9]{HJC17}.

Thus, SNR graphs of stabilizing Cox processes exhibit a supercritical phase if $r_{\mathrm B}=\ell^{-1}(\tau N_0)$ is large enough, and the critical intensity tends to zero as $r_{\mathrm B} \to \infty$. Hence, any intensity $\lambda>0$ is SNR-supercritical for $r_{\mathrm B}$ sufficiently large. That is, percolation can be obtained via reducing $N_0$ or $\tau$, or via increasing the signal power parameter $P$ that is fixed to the value $1$ in the present paper for simplicity. In practice, it depends on technological development and physical constraints whether such improvements are possible. We note that the paper \cite{BY13} worked under the assumption that $N_0$ an $\tau$ are fixed and thus formulated its results for large $r_{\mathrm B}$ in terms of large $P$.

If $d=2$, then\index{point process!Poisson} for the Poisson point process, Theorem~\ref{theorem-fromDF06} guarantees that\index{percolation!SINR!for Poisson point processes}\index{percolation!SINR!in two dimensions} $\lambda_{N_0,\tau}=\lambda_{\mathrm c}(r_{\mathrm B})<\infty$ for all $r_{\mathrm B}>0$. This relies on the\index{Russo--Seymour--Welsh type result|textbf} Russo--Seymour--Welsh type result \cite[Corollary 4.1]{MR96} that for $r>0$ and $\lambda>\lambda_{\rm c}(r)$,\index{crossing probabilities}\index{Boolean model!Poisson} $3n \times n$ rectangles are crossed by some cluster of the Poisson--Boolean model $X^\lambda \oplus B_{r/2}$ in the hard direction with probability tending to 1 as $n \to \infty$. This result is formulated more precisely and slightly more generally as follows (cf.~\cite[Theorem 2.7.1]{FM07}). For $n,\alpha>0$, let $\mathfrak R(\alpha,n)$ denote the rectangle $[0,\alpha n] \times [0, n] \subset \R^d$.

\begin{theorem}[MR96, FM07]\label{theorem-MR96}\index{crossing probabilities}
For $\lambda>0$, let $X^\lambda$ be a homogeneous Poisson point process on $\R^2$. Fix $r>0$. For $n>0$ and $\alpha>1$, let $C(\alpha,n)$ denote the event that $\mathfrak R(\alpha,n)$ is horizontally crossed by the Poisson--Boolean model $X^\lambda \oplus B_{r/2}$. That is, $C(\alpha,n)$ is defined as the event that there exists a connected component $\mathcal C$ of $X^\lambda \oplus B_{r/2}$ such that for both vertical sides $\lbrace 0 \rbrace \times [0,n]$ and $\lbrace \alpha n \rbrace \times [0,n]$ of $\mathfrak R(\alpha,n)$, there exists a point in $\mathcal C \cap X^\lambda$ having distance less than $r/2$ from that side. If $\lambda>\lambda_{\mathrm c}(r)$, then for any $\alpha>1$, we have that
$\lim_{n \to \infty} \mathbb P(C(\alpha,n))=1$.
\end{theorem}

Now, in the\index{limit!low density coupled with large radius} coupled limit $\lambda \downarrow 0, r \uparrow \infty, \lambda r^d=\varrho>0$, the rescaled\index{point process!Cox} Cox process $r^{-1} X^\lambda$ converges weakly to a\index{point process!Poisson} Poisson point process with intensity $\varrho$ \cite[Section~2.2.2]{HJC17}. Further, using arguments of \cite[Section~7.1]{HJC17}, we will see that for fixed $n$, the probability that the Boolean model of the Cox process crosses a $3nr \times nr$ rectangle in a given direction converges to the probability that the limiting Poisson--Boolean model crosses a $3n \times n$ rectangle in the same direction. These together with the stabilization of $\Lambda$ give us an opportunity to map the SINR graph to a\index{percolation process!renormalized} renormalized percolation process, using the construction of \cite[Section 3]{DF06} involving crossing probabilities. Moreover, if $\Lambda$ is stabilizing, then\index{interference control} interferences can be controlled similarly to the proof of Theorem~\ref{prop-firstpercolation} under the assumption \eqref{expmomentscond} on the stationary intensity measure $\Lambda$, using this renormalized percolation process. These imply that $\lambda_{N_0,\tau}<\infty$ if $r_{\mathrm B}$ is large. Actually, it is even true that any $\lambda>0$ exceeds $\lambda_{N_0,\tau}$ if $r_{\mathrm B}$ is sufficiently increased. Similarly, in higher dimensions $d \geq 3$, one can use discrete percolation arguments of \cite[Section 5.2]{HJC17} in order to verify an analogous assertion.



\begin{proposition}\label{prop-stabilizingsupercritical}\index{percolation!SINR!for Cox point processes}\index{percolation!SINR!for large connection radii}
Let $d \geq 2$ and $\lambda>0$, and let $\Lambda$ be stabilizing.
If $\supp(\ell)=[0,\infty)$ and assumption \eqref{expmomentscond} of Theorem~\ref{prop-firstpercolation} holds, then there exists $r_0 \geq 0$ such that if $r_{\mathrm B} \geq r_0$, then $\lambda_{N_0,\tau}<\lambda$.
\end{proposition}
We note that while $\lambda_{N_0,\tau}<\infty$ follows from the mere assumption that $r_{\mathrm B}$ is large, the function $\lambda \mapsto \gamma^*(\lambda)$ depends on finer details of the parameters $\lambda,\tau$, and $N_0$. E.g., for all $\lambda>0$, $\gamma^*(\lambda) \leq 1/\tau$ holds thanks to the degree bounds, see Section~\ref{sec-gammalambda}.

We will prove Proposition~\ref{prop-stabilizingsupercritical} in Section~\ref{sec-stabilizingsupercritical} and discuss its applicability to the main examples in Section~\ref{sec-examplesphasetransition}. 
Note that unlike Theorem~\ref{prop-firstpercolation}, Proposition~\ref{prop-stabilizingsupercritical} does not tell about the case when\index{path-loss function!compactly supported}\index{support!of the path-loss function} $\ell$ has compact support. Indeed, in that case, $r_{\mathrm B}$ cannot be increased arbitrarily, and it may happen that $\lambda_{\mathrm c}(r_{\mathrm B})=\infty$ for all $r_{\mathrm B}$ such that $\ell(r_{\mathrm B})>0$. Then, SINR graphs also do not percolate for any possible $r_{\mathrm B}<\sup\supp(\ell)$ and $\lambda>0$, $\gamma \geq 0$.

Although apart from the two-dimensional Poisson case we do not know whether $\lambda_{\mathrm c}(r_{\mathrm B})=\lambda_{N_0,\tau}$ holds for given values of the parameters, Proposition~\ref{prop-stabilizingsupercritical} implies at least that both critical intensities tend to zero as $r_{\mathrm B} \to \infty$. This assertion relies on the well-known \emph{scale invariance} of Poisson--Boolean models and Poisson--Gilbert graphs, cf.~\cite[Section~3]{DF06}.

\subsubsection{The case of no environmental noise}\label{sec-noIse}\index{noise}\index{percolation!SINR!for zero noise}
We now consider the case $N_0=0$. 
We fix $\tau>0$. Since for any $\tau,a>0$ and $\gamma>0$, one has\index{graph!SINR!decreasing in $N_0$} $g_{(\gamma,a,\tau)}(X^\lambda) \preceq g_{(\gamma,0,\tau)}(X^\lambda)$, it follows that \[ \lambda_{0,\tau} \leq \inf_{a>0} \lambda_{a,\tau}. \numberthis\label{trivialN0=0bound} \]


In the Poisson case for $d=2$, \cite[Section~3.4]{DF06} claimed that $\lambda_{0,\tau}=0$ and argued that this can be shown analogously to the statement of Theorem~\ref{theorem-fromDF06} that $\lambda_{N_0,\tau}<\infty$ for all $N_0>0$, and that the only difference is that there is no Boolean threshold. We now show that this claim is true if $\ell$ has unbounded support, but it fails in most of the relevant cases, in particular also in the two-dimensional Poisson case, if $\supp(\ell)$ is compact.

Let $\ell$ be such that $\supp~ \ell=[0,\infty)$. As for the case $d=2$ and $\Lambda \equiv |\cdot|$, let us fix $\tau>0$, and let $\lambda>0$ be arbitrary. By the\index{scale invariance!for Poisson--Boolean models} scale invariance of Poisson--Boolean models
and the fact that $\lambda_{\mathrm c}(1) \in (0,\infty)$, it follows that any $\lambda>0$ satisfies $\lambda>\lambda_{\mathrm c}(r)$ for all sufficiently large $r>0$. Choosing $r_{\mathrm B}(a)=\ell^{-1}(\tau a)$, we see that $r_{\mathrm B}(a)$ is well-defined for all sufficiently small noise powers $a>0$, and $r_{\mathrm B}(a) \to \infty$ as $a \downarrow 0$. The proof of \cite[Theorem~1]{DF06} implies that $\lambda_{a,\tau}=\lambda_{\mathrm c}(\ell^{-1}(\tau a))$ whenever the right-hand side of this equation is well-defined. Thus, $g_{(\gamma,a,\tau)}(X^\lambda)$ percolates almost surely for all $\gamma,a$ sufficiently small, and hence so does $g_{(\gamma,0,\tau)}(X^\lambda) \succeq g_{(\gamma,a,\tau)}(X^\lambda)$. 

Now, for $d \geq 2$, in the general Cox case, if $\supp(\ell)$ is unbounded, then letting $N_0 \downarrow 0$ is equivalent to letting $r_{\mathrm B} \to \infty$. Since $g_{(\gamma,N_0,\tau)}(X^\lambda) \preceq g_{(\gamma,0,\tau)}(X^\lambda)$ for any $N_0>0$, Proposition~\ref{prop-stabilizingsupercritical} implies that \eqref{trivialN0=0bound} holds and its right-hand side equals 0
if $\Lambda$ is stabilizing, $\supp(\ell)$ is unbounded, and \eqref{expmomentscond} holds. In contrast, if $\supp(\ell)$ is bounded, then for any $d \geq 2$, $\lambda_{0,\tau}=0$ is only true in the pathological case $\lambda_{\mathrm c}(r_{\max})=0$, in particular it never occurs if $\Lambda$ is stabilizing. 
\begin{corollary}\index{phase transition!in the SINR graph!without noise}
If $d \geq 2$ and $r_{\max}:=\sup\supp(\ell)$ is finite, then $\lambda_{0,\tau} \geq \lambda_{\mathrm c}(r_{\max})$.
\end{corollary}
The proof of the corollary uses an argument similar to the one in \cite[Section~3.4.2]{BY13}.
\begin{proof}
The statement is trivial if $\lambda_{\mathrm c}(r_{\max})=0$. Else, note that for any $\lambda>0$ and $X_i,X_j \in X^\lambda$, if $|X_i-X_j| \geq r_{\max}$, then $\SINR(X_i,~X_j,~X^\lambda)=0$. Hence, $g_{(\gamma,0,\tau)}(X^\lambda) \preceq g_{r_{\max}}(X^\lambda)$ for any $\gamma>0$. Choosing $0<\lambda < \lambda_{\mathrm c}(r_{\max})$, with probability 1, $g_{(\gamma,0,\tau)}(X^{\lambda})$ does not percolate for any $\gamma>0$.
\end{proof}
Thus, since $\Lambda \equiv |\cdot|$ is stabilizing, it follows that \cite[Corollary~1]{DF06} is false for all choices of $\ell$ with compact support.

\subsection{Estimates on the critical interference cancellation factor}\label{sec-gammalambda}
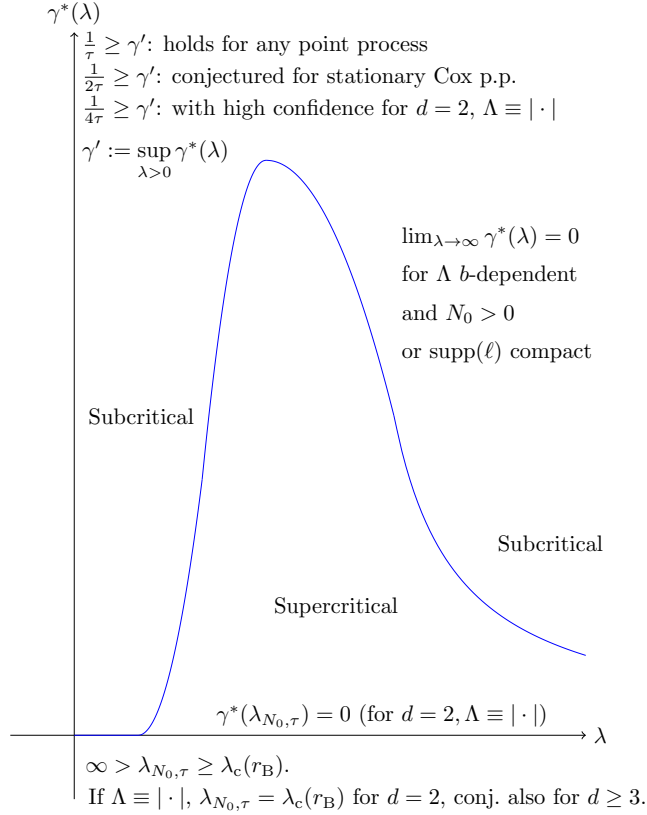
\begin{figure}
\scalebox{.85}{\begin{tikzpicture}[scale=1]
\draw[->] (-1,0) -- (8,0) node[right] {$\lambda$};
\draw[->] (0,-1) -- (0,11) node[above] {$\gamma^*(\lambda)$};
\draw[domain=0:1,smooth,variable=\x,blue] plot ({\x},{0});
\draw[domain=1:2,smooth,variable=\x,blue] plot ({\x},{4*\x^2-8*\x+4});
\draw[domain=2:3,smooth,variable=\x,blue] plot ({\x},{-5*\x^2+30*\x-36});
\draw[domain=3:5,smooth,variable=\x,blue] plot ({\x},{6*\x-\x^2});
\draw[domain=5:8,smooth,variable=\x,blue] plot ({\x},{5/(\x-4)});
\draw (0.05,-0.5) node[right] {$\infty>\lambda_{N_0,\tau} \geq \lambda_{\rm c}(r_{\rm B})$.};
\draw (0.1,-1) node[right] {If $\Lambda \equiv |\cdot|$,  $\lambda_{N_0,\tau} = \lambda_{\rm c}(r_{\rm B})$ for $d=2$, conj.~also for $d \geq 3$.};
\draw (3,2) node[right] {Supercritical};
\draw (2.1,0.3) node[right] {$\gamma^*(\lambda_{N_0,\tau})=0$ (for $d=2,\Lambda \equiv |\cdot|$)};
\draw (5,7.8) node[right] {$\lim_{\lambda \to \infty} \gamma^*(\lambda)=0$};
\draw (5,7.2) node[right] {for $\Lambda$ $b$-dependent};
\draw (5,6.6) node[right] {and $N_0>0$};
\draw (5,6) node[right] {or $\supp(\ell)$ compact};
\draw (6.5,3) node[right] {Subcritical};
\draw (0.1,5) node[right] {Subcritical};
\draw (0,9) node[right] {$\gamma':=\sup\limits_{\lambda>0} \gamma^*(\lambda)$};
\draw (0,10.3) node[right] {$\frac{1}{2\tau} \geq \gamma'$: conjectured for stationary Cox p.p.};
\draw (0,10.8) node[right] {$\frac{1}{\tau} \geq \gamma'$: holds for any point process};
\draw (0,9.8) node[right] {$\frac{1}{4 \tau} \geq \gamma'$: with high confidence for $d=2$, $\Lambda \equiv |\cdot|$};
\end{tikzpicture}}
\caption{Already proven (in black) and conjectured (in blue) properties of the $\lambda$--$\gamma^*(\lambda)$ phase diagram of the SINR graph of a Cox point process, in case $N_0,\gamma,\tau>0$ for $\Lambda$ asymptotically essentially connected under the condition \eqref{ellboundedcond} or \eqref{expmomentscond}, or for $r_{\rm B}$ sufficiently large for $\Lambda$ stabilizing under the condition \eqref{expmomentscond}. The question of the continuity of $\lambda \mapsto \gamma^*(\lambda)$ and the one of the uniqueness of its (local) maximum are open in general, and also whether its value at $\lambda_{N_0,\tau}$ equals zero. We note that in the two-dimensional Poisson case, according to \cite[Theorem 4.5]{MR96}, one has that $\P(g_{r_\mathrm B}(X^{\lambda_\mathrm c(r_\mathrm B)})\text{ percolates})=0$, and this together with Theorem~\ref{theorem-fromDF06} immediately implies that $\gamma^*(\lambda_{N_0,\tau})=0$.}\label{figure-gammastarlambda}
\end{figure}

\index{interference cancellation factor!estimates on the critical value}
In the Poisson case $\Lambda \equiv |\cdot|$ for $d=2$, \cite{DBT03, FM07} derived the following bounds on the critical interference cancellation factor $\gamma^*(\lambda)$ defined in \eqref{gammastarlambda}.
\begin{enumerate}[(A)]
\item\label{first-Poissonbound} $\forall \lambda>0$, $\gamma^*(\lambda) \leq \frac{1}{\tau}$.
\item\label{second-Poissonbound} $\gamma^*(\lambda) = \Ocal(1/\lambda)$ as $\lambda \to \infty$.
\item\label{third-Poissonbound} If $\ell$ has bounded support, then $\gamma^*(\lambda) = \Omega(1/\lambda)$ as $\lambda \to \infty$.
\end{enumerate}
\eqref{first-Poissonbound} implies that $\lambda \mapsto \gamma^*(\lambda)$ is bounded.
In Section~\ref{sec-generalgammabound} we recover this bound for any stationary Cox point process and present conjectures regarding its possible improvements. 
For the Cox case, in Section~\ref{sec-lambdadependentgammabound} we provide sufficient conditions under which \eqref{second-Poissonbound} holds or at least $\gamma^*(\lambda)$ tends to 0 as $\lambda \to \infty$, while in Section~\ref{sec-Omega(1/lambda)discussion} we investigate generalizations of \eqref{third-Poissonbound}. Figure~\ref{figure-gammastarlambda} visualizes our results and conjectures.


\subsubsection{Intensity-independent bounds}\label{sec-generalgammabound}\index{interference cancellation factor!estimates on the critical value!uniform upper bound}
In the Poisson case, \eqref{first-Poissonbound} is a consequence of the fact\index{phase transition!in the SINR graph!in view of degree bounds} \cite[Theorem~1]{DBT03} that SINR graphs with $\gamma>0$ have bounded degrees. This assertion generalizes to any dimension and any simple point process \cite[Section 4.1.4.1]{T19}.


\begin{proposition}\label{prop-gammalambdaestimate1}\index{degree bounds!for general simple point processes}
Let $\lambda>0$, $\tau>0$, and $N_0>0$. Then, almost surely,\index{degree bounds|textbf}
\[\forall \gamma>0, \forall i \in \N, \quad X_i\text{ has in-degree less than } 1+\smfrac{1}{\tau\gamma}\text{ in }g^{\rightarrow}_{(\gamma,N_0,\tau)}(X^\lambda). \numberthis\label{directeddegreebound} \] 
In particular, for all $\lambda>0$, $\gamma^*(\lambda) \leq \frac{1}{\tau}$, and $\P(g_{(\frac{1}{\tau},N_0,\tau)}(X^\lambda)\text{ percolates})=0$.
\end{proposition} For $N_0=0$, the same proof implies the same assertion apart from the non-percolation for $\gamma=1/\tau$. 
The proof of the bound \eqref{directeddegreebound} is analogous to the one of \cite[Theorem~1]{DBT03}. We note that it even holds in one dimension, and among the properties of $\ell$ it only uses that $\ell(|X_i-X_j|)>0$ holds if there is an edge from $X_i$ to $X_j$ or from $X_j$ to $X_i$ in the directed SINR graph. The
arguments of its proof can also be used in order to derive stronger degree bounds if $N_0>0$ and to show that also the out-degrees in $g^{\rightarrow}_{(\gamma,N_0,\tau)}(X^\lambda)$ are bounded if $\ell$ has unbounded support; we refrain from presenting here the details.

By \eqref{directeddegreebound}, if $\gamma \geq \smfrac{1}{\tau}$,\index{degree bounds!equal to one} degrees  in $g_{(\gamma,N_0,\tau)}(X^\lambda)$  are at most 1, and thus all\index{cluster!in the SINR graph} clusters of $g_{(N_0, \gamma,\tau)}(X^\lambda)$ are pairs or isolated points. This implies lack of percolation. We also expect that there is no infinite cluster if $\gamma \in [\smfrac{1}{2\tau},\smfrac{1}{\tau})$, where\index{degree bounds!equal to two} the degree bound is 2, for a large class of point processing including the stationary Poisson one. Indeed, in this regime, all clusters are isolated points, finite cycles or (possible in one or two directions infinite) paths. This reminds of one-dimensional percolation models, which are very often subcritical (see e.g. \cite[Section 3.2]{MR96}).

The degree constraints also relate SINR graphs to certain\index{graph!$k$-nearest neighbour!bidirectional|textbf} $k$-nearest neighbour graphs. For $k \in \N$ and $\lambda>0$, the \emph{bidirectional $k$-nearest neighbour graph} $g_{\mathbf B}(k,X^\lambda)$ is defined as the undirected graph where $X_i,X_j \in X^\lambda$, $i \neq j$, are connected by an edge if and only if they are mutually among the $k$ nearest neighbours of each other. It is easy to see that this graph is almost surely well-defined for stationary Cox point processes. We have the following lemma, which can also be generalized for a large class of simple point processes, cf.~\cite[Section 4.1.4.1]{T19}. 
\begin{lemma}\label{lemma-neighbourrelations}\index{graph!SINR!relation to $k$-nearest neighbour graphs}
For any stationary $\Lambda$, for any $\lambda>0$ and $k \in \N$, if $\tau,\gamma>0$ and $N_0 \geq 0$ are such that almost surely, all in-degrees in $g^{\rightarrow}_{(\gamma,N_0,\tau)}(X^\lambda)$ are at most $k \in \N$, then  $g_{(\gamma,N_0,\tau)}(X^\lambda) \preceq g_{\mathbf B}(k,X^\lambda)$.
\end{lemma}
The proof of Lemma~\ref{lemma-neighbourrelations} is immediate, therefore we omit it. We use this lemma in order to derive a conjecture for\index{point process!Poisson} the two-dimensional Poisson case. In this case \cite{BB08} studied the graph $g_{\mathbf B}(k,X^1)$, which has the same distribution as $\lambda^{1/2} g_{\mathbf B}(k,X^\lambda)$ for all $\lambda>0$. In particular, $\mathbb P(g_{\mathbf B}(k,X^1)\text{ percolates})=\mathbb P(g_{\mathbf B}(k,X^\lambda)\text{ percolates})$ for all $\lambda>0$.
By \cite[Section~3]{BB08}, \index{high-confidence result}with high confidence, $g_{\mathbf B}(k,X^1)$ percolates only if $k \geq 5$. That is, this assertion follows once one proves that certain high-dimensional integrals exceed certain deterministic values, but so far the integrals have only been evaluated using Monte Carlo methods. This is more than simulations but less than a proof. If the result holds, then by \eqref{directeddegreebound} and Lemma~\ref{lemma-neighbourrelations}, it implies the following improvement of \cite[Theorem~1]{DBT03}.
\begin{conjecture}
Let $\Lambda \equiv |\cdot|$ and $d=2$. Then for any $N_0 \geq 0$ and $\lambda>0$,
$\gamma^*(\lambda) \leq \smfrac{1}{4\tau}$,
and $\P(g_{(\frac{1}{4\tau},N_0,\tau)}(X^\lambda)\text{ percolates})=0$.
\end{conjecture}
\begin{figure}
\includegraphics[scale=0.20]{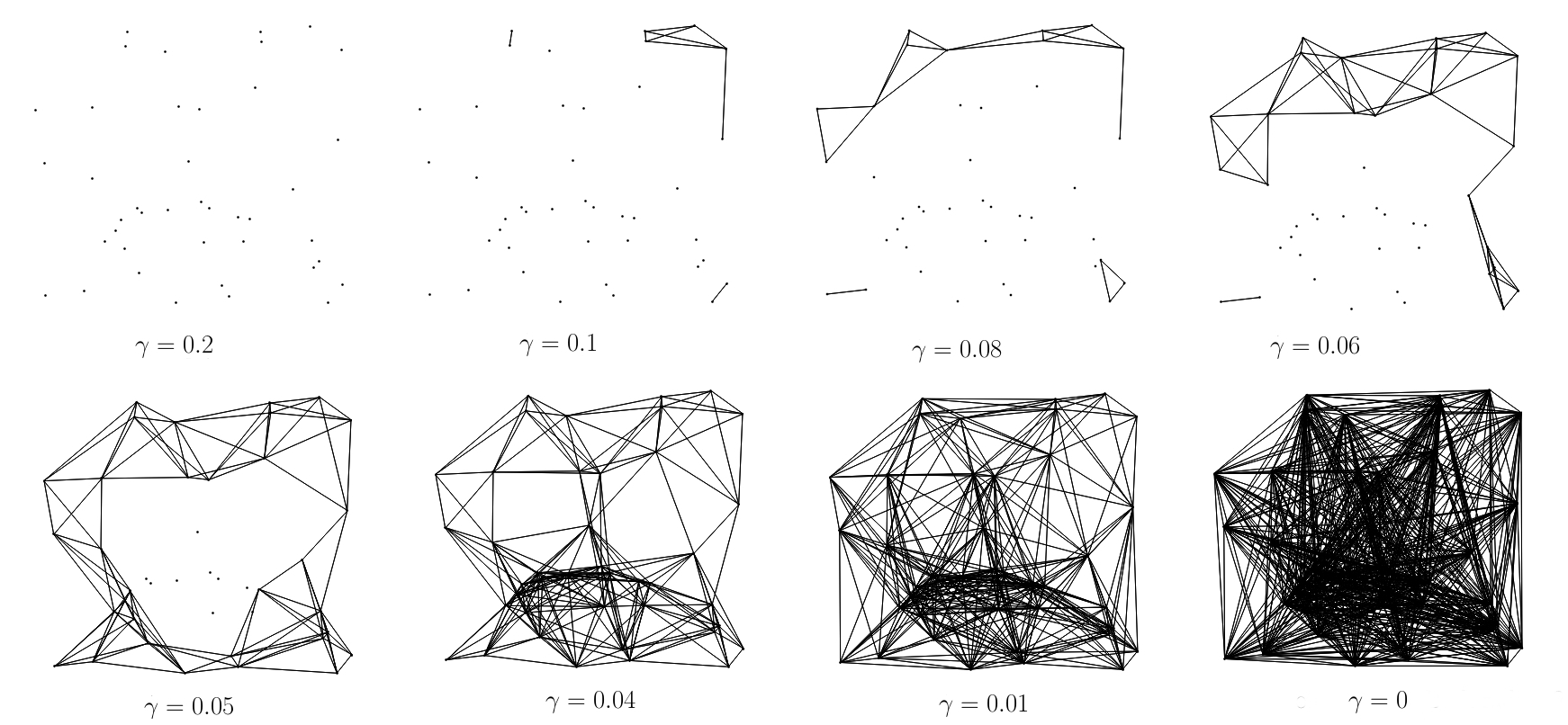}
\caption{$\SINR$ graphs for $d=2$ and $\Lambda \equiv |\cdot|$ restricted to $[0,1]^2$ with different values of $\gamma$, where $N_0=2,\tau=1=1,\ell(r)=\min\lbrace 100,r^{-4} \rbrace$, $\lambda=40$. The realization for $\gamma=\frac{1}{5\tau}$ has no edges, and the one for $\gamma=\frac{1}{10\tau}$ is still highly disconnected. The ones for $\gamma \leq \frac{1}{25\tau}$ are connected, but the effect of bounded degrees is still prominent for $\gamma=\frac{1}{100\tau}$ in comparison with the almost complete graph corresponding to $\gamma=0$. }\label{figure-gammas}
\end{figure}
\begin{figure}
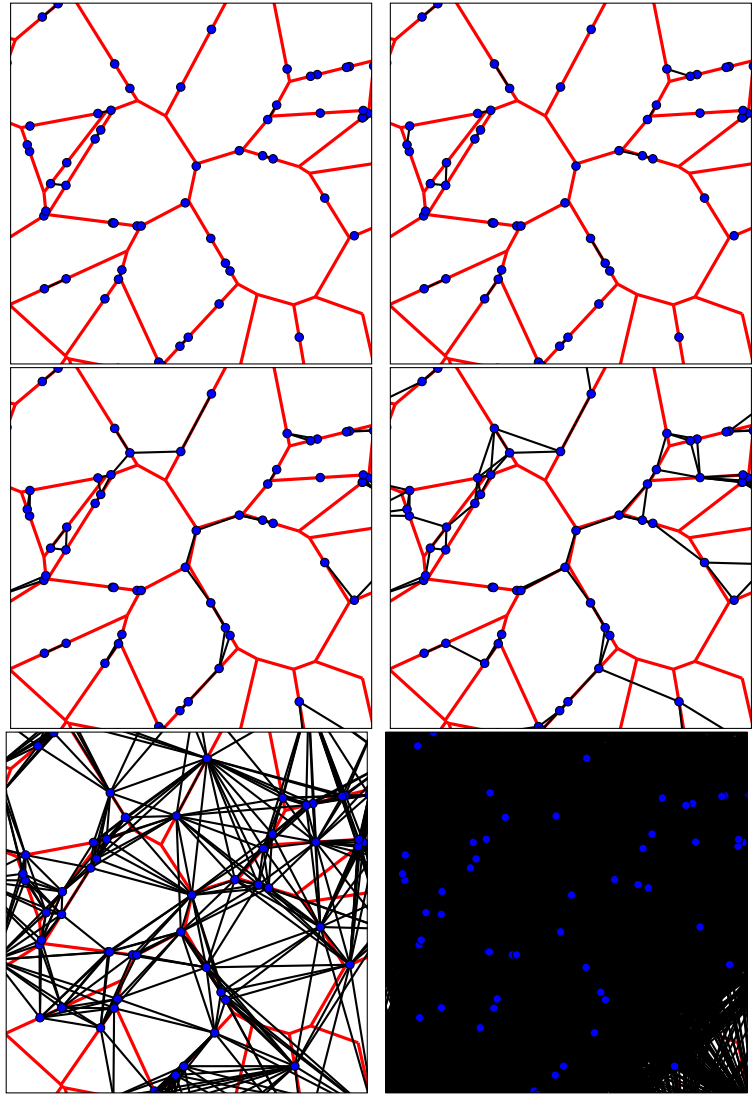

\input{gamma=1.tex}
\input{gamma=0.5.tex}
\input{gamma=0.15.tex}
\input{gamma=0.07.tex}
\input{gamma=0.01.tex}
\input{gamma=0.tex}
\caption{$\SINR$ graphs for $d=2$ and $\Lambda$ given by a Poisson--Voronoi tessellation restricted to $[-0.1,1.1]^2$, from which we only view the part corresponding to the users in $[0,1]^2$ in order to handle boundary effects better. The parameters are $N_0=\tau=1,\ell(r)=\min \lbrace \eps, r^{-3} \rbrace$ for a certain $\eps$ that is smaller than the nearest neighbour distance of the given realization of the Cox process, and $\lambda$ much larger than $\lambda_{\rm c}(r_{\rm B})$.
(1) $\gamma=\frac{1}{\tau}$: the degree bound is 1, the SINR graph consists of pairs and isolated points. Note that the value of $\gamma$ where the last edge disappears is between $\frac{4\times 10^5}{\tau}$ and $\frac{5\times 10^5}{ \tau}$.
(2) $\gamma=\frac{1}{2\tau}$: the degree bound is 2, the SINR graph contains some points with degree 2 but hardly any cycles of length larger than 2, and it is highly disconnected.
(3) $\gamma=\frac{15}{100\tau}$: the degree bound is 7. While the graph is still clearly subcritical, the first larger cycles have already arised.
(4) $\gamma=\frac{7}{100\tau}$ seems to be close to the critical value $\gamma^*(\lambda)$: most of the realization consists of the two biggest clusters, which do not yet connect up in $[0,1]^2$. The degree bound is 15, the average degree is about 3.18. (5) $\gamma=\frac{1}{100\tau}$: the degree bound is 100, the number of points in the realization is 119. The graph clearly appears to be supercritical but is still much sparser than for $\gamma=0$. (6) $\gamma=0$. The SNR graph is almost complete, it contains 6945 edges out of the 7021 ones of the complete graph.}
\label{figure-PVTgammas}
\end{figure}

Simulations suggest that the maximum of $\lambda \mapsto \gamma^*(\lambda)$ is even lower than $\smfrac{1}{4\tau}$, cf.~Figures~\ref{figure-gammas} and \ref{figure-PVTgammas}. Conversely, Theorem~\ref{prop-firstpercolation}\eqref{firstsupercritical} implies that for $d \geq 2$, $g_{\mathbf B}(k,X^1)$ percolates for all $k$ sufficiently large. This was proven in \cite{BB08} for $d=2$ and $k \geq 15$, and it is intuitively quite clear that this implies the same statement for any $d \geq 3$ for $k$ sufficiently large, although this was not explicitly stated in \cite{BB08}.

\subsubsection{Upper bounds for large intensities}\label{sec-lambdadependentgammabound}\index{interference cancellation factor!estimates on the critical value!upper bound for large $\lambda$}\index{percolation!SINR!in two dimensions}
For\index{b-dependence@$b$-dependence} $b$-dependent Cox processes in $d \geq 3$ dimensions, we recover \eqref{second-Poissonbound} in a weaker form. 
Namely, any $\gamma>0$ becomes subcritical for large $\lambda$ whenever the SINR graph has\index{graph!SINR!with bounded edge length} bounded edge length. 
\begin{proposition}\label{prop-gammalambda0}
If $\Lambda$ is $b$-dependent, $N_0 \geq 0$, $\tau>0$, further, $N_0>0$ or $\ell$ has bounded support, then 
\[ \lim_{\lambda \to \infty} \gamma^*(\lambda) = 0. \numberthis\label{largelambdasmallgamma}\]
\end{proposition}
We will prove Proposition~\ref{prop-gammalambda0} in Section~\ref{sec-gammalambda0proof}. 

Further, for $d=2$, \eqref{second-Poissonbound} stays true for $b$-dependent Cox processes for which $\Lambda(Q_{\delta})$ is bounded away from 0 for small enough $\delta>0$. 
\begin{proposition}\label{prop-bdependentOlambda}\index{measure!intensity!locally bounded away from zero}
If $d=2$, $N_0,\tau>0$, $\Lambda$ is $b$-dependent, and $\Lambda(Q_{\delta/2})$ is bounded away from 0 for some $\delta>0$ such that $\ell(\delta)>\tau N_0$, then  
as $\lambda \to \infty$,  
$\gamma^*(\lambda) = O (1/\lambda)$. 
\end{proposition}
The proof of Proposition~\ref{prop-bdependentOlambda} will be carried out in Section~\ref{sec-bdependentOlambdaproof}. The reason why this proposition is restricted to $d=2$ is that its proof uses that in a certain $b$-dependent site percolation model, the origin not being contained in an infinite cluster is equivalent to the origin being surrounded by a circuit of closed sites, which has no clear analogue for higher dimensions. The applicability of the results of this section to the main examples will be discussed in Section~\ref{sec-examplesgammalambda}.


\subsubsection{Lower bounds for large intensities}\label{sec-Omega(1/lambda)discussion}\index{interference cancellation factor!estimates on the critical value!lower bound for large $\lambda$}\index{percolation!SINR!in two dimensions}
In \cite[Section~III-C]{DBT03}, \eqref{third-Poissonbound} was verified for the Poisson case for $d=2$ and compactly supported $\ell$. It can easily be generalized to a class of $b$-dependent Cox point processes as follows.
\begin{corollary}\label{cor-Omega1/lambda}
Let $d=2$, $\supp(\ell)$ bounded, and let $\Lambda$ be $b$-dependent such that\index{measure!intensity!locally bounded away from zero} $\Lambda(Q_{\eta})$ is bounded away from 0 for some $\eta>0$. Then we have $\gamma^*(\lambda)=\Omega(1/\lambda)$ as $\lambda \to \infty$.
\end{corollary}
The proof of Corollary~\ref{cor-Omega1/lambda} will be sketched in Section~\ref{sec-Omega(1/lambda)proof}. This follows the lines of the original proof of \cite[Theorem 2]{DBT03}, using some additional observations. Since this proof involves a dual lattice argument that is not applicable in higher dimensions, the corollary remains restricted to $d=2$. In Section~\ref{sec-examplesgammalambda} we will discuss the applicability of Corollary~\ref{cor-Omega1/lambda} to the main examples.


\subsection{Applicability of the results to the main examples}\label{sec-applicability}
\subsubsection{Phase transitions}\label{sec-examplesphasetransition}
We now consider each of the relevant examples of $\Lambda$ from \cite{HJC17} recalled in Section~\ref{sec-Gilbertmodeldef} and discuss the applicability of Theorem~\ref{prop-firstpercolation} and Proposition~\ref{prop-stabilizingsupercritical} to them. For the sake of brevity, we will tacitly assume that $N_0>0$. The case $N_0=0$ can be handled according to Section~\ref{sec-noIse}.

Let us note that in the case of a modulated Poisson point process, $\Lambda$ is\index{b-dependence@$b$-dependence!examples} $b$-dependent for some $b>0$ also if the stationary random set $\Xi$ is a generalized Boolean model with bounded grains (cf.~\cite[Chapter 3]{BB09}). That is, $\Xi=\bigcup_{i \in \N} \Gcal(X_i)$, where the grains $\Gcal(X_i)$ (which depend on the position of the Cox point $X_i$ but not on the other points of $X^\lambda$) are such that there exists a compact set $K \subseteq \R^d$ such that $\Gcal(X_i)-X_i \subseteq K$ for all $\R^d$. The classical Boolean model corresponds to the special case $\Gcal(X_i)=B_{r/2}(X_i)$.

Now, all examples are stabilizing and therefore they exhibit a subcritical phase by Theorem~\ref{prop-firstpercolation}\eqref{firstsubcritical}, apart from general modulated Poisson point processes where $\Xi$ is not $b$-dependent.

For the modulated Poisson point process with $\Xi$ being a generalized Boolean model with bounded grains, $\Lambda$ is $b$-dependent and $\Lambda(Q_1)$ is bounded. Further, if $\lambda_1,\lambda_2>0$, or if $\lambda_1>\lambda_2=0$ and $\Xi$ is a supercritical Boolean model, then $\Lambda$ is\index{asymptotic essential connectedness!examples} asymptotically essentially connected. In these cases, $\lambda_{N_0,\tau}<\infty$ holds for any $r_{\mathrm B}>0$ under the general Assumption~($\ell$) on $\ell$. In particular, this covers the Poisson case $\Lambda \equiv |\cdot|$. Further, by stabilization, $\lambda_{N_0,\tau}<\infty$ holds for large $r_{\mathrm B}$ also if either $\lambda_1$ or $\lambda_2$ is zero, in case $\ell$ has unbounded support and satisfies Assumption~($\ell$). It is easy to see that if $\Xi$ is a Poisson--Boolean model and $\lambda_1>\lambda_2=0$,\index{phase!supercritical!examples for absence for stabilizing $\Lambda$ and small $r$}\index{stabilization!examples!showing lack of supercritical phase for small $r$}\index{phase!supercritical!examples for absence for stabilizing $\Lambda$ for small $r$} then there are cases where $\lambda_{\mathrm c}(r_{\mathrm B})=\infty$ holds for small $r_{\mathrm B}>0$, cf.~Figure~\ref{figure-MPP}. Indeed,\index{point process!Poisson, modulated!examples for lack of supercritical phase} if the Poisson--Boolean model is subcritical, then one can choose $r_{\mathrm B}$ so small that the\index{Boolean model!Cox} Cox--Boolean model $X^\lambda \oplus B_{r_{\mathrm B}/2}$ is still contained in a subcritical Poisson--Boolean model for any $\lambda>0$. Also for $\lambda_2>\lambda_1=0$, a supercritical phase may be missing. Indeed, e.g.~for $d = 2$ and $\lambda>0$, for any supercritical Poisson--Boolean model $\mathcal B(\lambda_0,r_0)$ with intensity $\lambda_0>0$ and radius $r_0>0$, there exists $r_1<r_0$ such that $\mathcal B(\lambda_0,r_1)$ has no unbounded vacant component \cite[Section~4.6]{MR96}. Then for $\Xi=\mathcal B(\lambda_0,r_0)$ and $\lambda>0$, let the Cox process $X^\lambda$ have intensity measure $\lambda\Lambda$, with $\Lambda=\lambda_2 \ind{1}_{\Xi^{\mathrm c}}|\cdot|$ satisfying $\E[\Lambda(Q_1)]=1$. Then for $r_{\rm B}>0$ small, for all $\lambda>0$, the Cox--Boolean model $X^\lambda \oplus B_{r_{\mathrm B}/2} $ is included in $\mathcal B(\lambda_0,r_1)^{\mathrm c}$ and thus has no unbounded cluster. 

For a general, not $b$-dependent $\Xi$, neither Theorem~\ref{prop-firstpercolation} nor Proposition~\ref{prop-stabilizingsupercritical} is applicable due to the possible lack of stabilization. However, $\Lambda(Q_1)$ is still bounded and $\Lambda$ is absolutely continuous, and therefore a subcritical phase exists for $\lambda_1,\lambda_2 \geq 0$ thanks to a comparison to a Poisson--Gilbert graph. Further, if $\lambda_1,\lambda_2>0$, then a similar comparison yields that $\lambda_{N_0,\tau}<\infty$ holds for any $r_{\mathrm B}$. These assertions were proven in \cite[Section 4.2.3.2]{T19}. 

For the\index{shot-noise field!examples for lack of supercritical phase} shot-noise field, let us recall that we only consider the case when the kernel $k$ is compactly supported and hence $\Lambda$ is $b$-dependent. For this intensity measure, it may again happen that $\lambda_{\mathrm c}(r_{\mathrm B})=\infty$ for small $r_{\mathrm B}>0$. Indeed, if the underlying Poisson point process $X_{\mathbf S}$ is such that its Boolean model  with connection radius $r/2=\diam~\supp~k$  is subcritical and also $r_{\mathrm B}$ is small, then the Cox--Boolean model $X^\lambda \oplus B_{r_{\mathrm B}/2}$ is included in a subcritical Poisson--Boolean model for any $\lambda>0$. Nevertheless, for any shot-noise field, $\Lambda$ is $b$-dependent and, although $\Lambda(Q_1)$ is not bounded, it has\index{exponential moments of the intensity!shot-noise field}\index{shot-noise field!existence of all exponential moments} all exponential moments thanks to\index{Campbell's theorem} Campbell's theorem \cite[Section~3.2]{K93}. Hence, for $\ell$ with unbounded support satisfying Assumption $(\ell)$, $\lambda_{N_0,\tau}<\infty$ holds for the shot-noise field if $r_{\mathrm B}$ is large.

On the other hand, there exist shot-noise fields that are asymptotically essentially connected. For example, take a supercritical Poisson--Boolean model $\Xi$, and denote its connection radius by $r$. Let now $k \colon \R^d \to [0,\infty)$ be defined as $k(x)=\mathds \eps_1 \{ |x|<r/2 \}$, where $\eps_1>0$ is uniquely chosen so that $\E[\Lambda(Q_1)]=1$ for the corresponding shot-noise field intensity $\Lambda$. Since $\supp(\Lambda)$ equals a supercritical Poisson--Boolean model, $\Lambda$ is asymptotically essentially connected.

Poisson--Voronoi and Poisson--Delaunay tessellations are asymptotically essentially connected, and thus by \eqref{ellboundedcond}, $\lambda_{N_0,\tau}<\infty$ holds for any $r_{\mathrm B}>0$ if $\ell$ has bounded support. They are neither $b$-dependent nor bounded, hence the question of existence of a supercritical phase for $\supp(\ell)$ unbounded remains open. On the other hand, thanks to the results of \cite{JT19}, for $d=2$ they satisfy the exponential moment condition in \eqref{expmomentscond}. 

\subsubsection{Estimates on the critical interference cancellation factor}\label{sec-examplesgammalambda}
Let us now discuss the applicability of Propositions~\ref{prop-gammalambda0} and \ref{prop-bdependentOlambda}, and Corollary~\ref{cor-Omega1/lambda} to the main examples. Each of them requires $b$-dependence, and therefore they are only applicable to the Poisson point process modulated by a generalized Boolean model with bounded grains and to the shot-noise field. For these two examples, Proposition~\ref{prop-gammalambda0} immediately applies. For $d=2$, Proposition~\ref{prop-bdependentOlambda} and Corollary~\ref{cor-Omega1/lambda} require also that $\Lambda(Q_\delta)$ be bounded away from 0 for some $\delta>0$, which only applies for the modulated Poisson point process with a $b$-dependent $\Xi$ and with $\lambda_1,\lambda_2>0$ (for which it holds for all $\delta>0$). 

\section{Proof and discussion of phase transitions}\label{sec-phasetransitionproofs}
This section includes the proofs of the results of Section~\ref{sec-phasetransitiondescription}. In particular, in Section~\ref{sec-firstpercolationproof} we verify Theorem~\ref{prop-firstpercolation}\eqref{firstsupercritical}. In Section~\ref{sec-additionalchallenges}, we comment on this proof, in particular on the interference control argument. 
Further, Section~\ref{sec-stabilizingproof} contains the proof of the results of Section~\ref{sec-stabilizing}:  in Section~\ref{sec-largersupercritical} we show how Corollary~\ref{cor-largersupercritical} can be derived from the results of \cite{HJC17}, whereas in Section~\ref{sec-stabilizingsupercritical}, using arguments of Section~\ref{sec-firstpercolationproof}, we verify Proposition~\ref{prop-stabilizingsupercritical}. 

\subsection{Proof and discussion of Theorem~\ref{prop-firstpercolation}}
\subsubsection{Proof of Theorem~\ref{prop-firstpercolation} (existence of supercritical phase)}\label{sec-firstpercolationproof}

For the proof we fix $N_0,\tau>0$. Now, for $\gamma \geq 0$ and $\lambda>0$, we use the simplified notation $g_{(\gamma)}(X^\lambda)=g_{(\gamma,N_0,\tau)}(X^\lambda)$ (until the end of the present section). Further, we assume that $\Lambda$ is asymptotically essentially connected. Thus, by \cite[Theorems 2.4, 2.6]{HJC17}, $\lambda_{\mathrm c}(r) \in (0,\infty)$ holds for all $r>0$. We recall that $g_{(0)}(X^\lambda)=g_{r_{\mathrm B}}(X^\lambda)$, cf.~\eqref{criticalradius}. The proof follows the four-step strategy that was outlined in Section~\ref{sec-positivenoise}. 


\begin{step}\label{step-latticemapping}
Mapping to a lattice percolation problem.
\end{step}
Let $r \in (\upsilon_0,r_{\mathrm B})$, such $r$ exists by \eqref{criticalradius} and \eqref{first-pathloss} -- \eqref{third-pathloss} in Assumption~($\ell$). Following \cite[Section~5.2]{HJC17}, for $n \geq 1$, we let $\Xi_n(nz)$ denote the union of all connected components of
$\supp(\Lambda_{Q_n(nz)} )$ that are of diameter at least $n/3$, and we say that a site $z \in \Z^d$ is $n$\emph{-good} if\index{percolation process!discrete!$b$-dependent}\index{percolation process!discrete!site}
\begin{enumerate}
\item\label{first-ngood1} $R(Q_n(nz))<n/2$,
\item\label{second-ngood1} $X^\lambda \cap \Xi_n(n z) \neq \emptyset$, and
\item\label{third-ngood1} for every $z' \in \Z^d$ with $|z-z'|_{\infty} \leq 1$ it holds that every  $X_i \in X^\lambda \cap \Xi_n(nz)$ and $X_j \in X^\lambda \cap \Xi_n(nz')$ are connected by a path in $g_r(X^\lambda) \cap Q_{6n}(nz)$.
\end{enumerate}
A site $z \in \Z^d$ is \emph{$n$-bad} if $z$ is not $n$-good. 

Next, for $a \geq 0$, we define a\index{path-loss function!shifted} ``shifted'' version $\ell_a$ of the path-loss function $\ell$, similarly to \cite{DF06}, which will be used in order to estimate interference values from above. Note that any point of $Q_{a}(x)$ is at distance at most $\smfrac{a\sqrt d}{2}$ away from the centre $x$ of $Q_{a}(x)$. We define $\ell_{a} \colon [0,\infty) \to [0,\infty)$ as follows
\[ \ell_{a}(r) = \ell(0) \ind{1}\Big\lbrace r< \frac{a\sqrt d}{2} \Big\rbrace + \ell\Big( r- \frac{a\sqrt d}{2}  \Big) \ind{1}\Big\lbrace r \geq  \frac{a\sqrt d}{2} \Big\rbrace. \numberthis\label{elladef} \]
Note that $\ell_0=\ell$. 
Now, we define the\index{shot-noise process} shot-noise processes
\begin{eqnarray*}
& I_{a}(x)=\sum_{X_i \in X^\lambda} \ell_{a}(|x-X_i|), \qquad I(x)=\sum_{X_i \in X^\lambda} \ell(|x-X_i|), \qquad x \in \R^d.
\end{eqnarray*}
Then $I_{0}(x)=I(x)$. By the triangle inequality, for $a \geq 0$, $I(x) \leq I_{a}(z)$ holds for any $z \in \Z^d$ and $x \in Q_a(z)$. 
Now, for $z \in \Z^d$, $n \geq 1$, and $M>0$, we define the following events 
\[
A_n(z) =\lbrace z\text{ is }n\text{-good}\rbrace,~~ B_{n,M}(z) =\lbrace I_{6n}(nz) \leq M \rbrace ,~~
 C_{n,M}(z)=A_n(z) \cap B_{n,M}(z). 
\]
\begin{step}\label{step-latticepercolation}
Percolation in the lattice. 
\end{step}
If $\lambda>0$ is sufficiently large, then for all $n,M$ sufficiently large, the process of $n$-good sites $z \in \Z^d$ such that $I_{6n}(nz) \leq M$ percolates with probability one (where $\Z^d$ is equipped with its nearest neighbour edges). This immediately follows by a\index{Peierls argument} Peierls argument (cf.~\cite[Section 1.4]{G99}) once we have verified that the following holds.
\begin{proposition}\label{prop-Cunionbound}
Under the assumption \eqref{ellboundedcond} or \eqref{expmomentscond} in Theorem~\ref{prop-firstpercolation}, for all sufficiently large  $n \geq 1$,  $\lambda=\lambda(n)>0$, and $M=M(\lambda,n)>0$, there exists a constant $q_C=q_C(\lambda,n,M) <1$ such that for any $L \in \N$ and pairwise distinct sites $z_1,\ldots,z_L\in \Z^d$, we have
\[ \mathbb P(C_{n,M}(z_1)^{\rm c} \cap \ldots \cap C_{n,M}(z_L)^{\rm c}) \leq q_C^L. \numberthis\label{qCexpdecay} \]
Moreover, for any $\eps>0$, we can choose $\lambda$, $n$, and $M$ large enough  such  that $q_C \leq \eps$. 
\end{proposition}
In order to verify this proposition we start with the results of \cite{HJC17} about the $n$-good sites.
\begin{lemma}\label{lemma-mygoodness}
For all sufficiently large $n \geq 1$ and $\lambda=\lambda(n) >0$, there exists $q_A=q_A(n,\lambda)<1$ such that for any $L \in\N$ and pairwise distinct sites $z_1,\ldots,z_L\in \Z^d$,
\[ \mathbb P(A_n(z_1)^{\rm c} \cap \ldots \cap A_n(z_L)^{\rm c})  \leq q_A^L. \numberthis\label{qAexpdecay} \]
Moreover, for any $\eps>0$ and for sufficiently large $n$, one can choose $\lambda$ so large that $q_A \leq \eps$.
\end{lemma}
\begin{proof}
In \cite[Section~5.2]{HJC17} it was shown that for\index{asymptotic essential connectedness}\index{asymptotic essential connectedness!implies existence of supercritical phase} asymptotically essentially connected $\Lambda$, the process of $n$-good sites is 7-dependent
. Moreover, for $z \in \Z^d$, we have
\[ \lim_{n \to \infty} \lim_{\lambda \to \infty} \mathbb P(A_n(z)^{\rm c})=0, \numberthis\label{goodnesswins} \]
where the convergence is uniform in $z \in \Z^d$. Let now $L \in \N$ and $z_1,\ldots,z_L \in \Z^d$  pairwise distinct. Let us write $[L]=\lbrace 1,\ldots,L\rbrace$. By 7-dependence, there exists $m \geq 1$ and a subset $\lbrace k_j \colon j=1,\ldots,m \rbrace$ of $[L]$ such that $A_n(z_{k_1}),\ldots,A_n(z_{k_m})$ are independent and $m \geq \frac{L}{8^d}$. Now, let $q'_A(n)=\limsup_{\lambda \to \infty} \mathbb P(A_n(o)^{\rm c})^{\frac{1}{8^d}}$. By \eqref{goodnesswins}, $q'_A(n)$ tends to zero as $n \to \infty$. Hence, for $n \geq 1$ sufficiently large, there exists $\lambda=\lambda(n)>0$ such that 
\[ \mathbb P(A_n(z_1)^{\rm c} \cap \ldots \cap A_n(z_L)^{\rm c} ) \leq  \mathbb P(A_n(z_{k_1})^{\rm c} \cap \ldots \cap A_n(z_{k_m})^{\rm c})  \leq \mathbb P(A_n(o)^{\rm c})^{\frac{L}{8^d}}\leq q_A^L, \]
where $q_A=2q'_A(n)$. This finishes the proof of the lemma.
\end{proof}

The main step of the proof of Proposition~\ref{prop-Cunionbound} is to prove the following assertion, in other words, to control the interferences.
\begin{proposition}\label{prop-qBexpdecay}\index{interference!control}
Under the assumption \eqref{ellboundedcond} or \eqref{expmomentscond} in Theorem~\ref{prop-firstpercolation}, for all sufficiently large $n \geq 1$, for all $\lambda>0$, and for all sufficiently large $M=M(n,\lambda)>0$, there exists a constant $q_B=q_B(n,\lambda,M) <1$ such that for any $L \in \N$ and pairwise distinct sites $z_1,\ldots,z_L\in \Z^d$, we have
\[ \mathbb P(B_{n,M}(z_1)^{\rm c} \cap \ldots \cap B_{n,M}(z_L)^{\rm c}) \leq q_B^L. \numberthis\label{qBexpdecay} \]
Moreover, for any $\eps>0$, for all large enough $n \geq 1$ and for all $\lambda>0$, we can choose $M$ large enough  such  that $q_B \leq \eps$.
\end{proposition}
This proposition is formally analogous to \cite[Proposition 2]{DF06} (apart from the additional technical condition that $n$ has to be large enough). 
The proof of Proposition~\ref{prop-qBexpdecay} is however more involved; it is postponed until Step~\ref{step-interferencecontrol}. Given Lemma~\ref{lemma-mygoodness} and Proposition~\ref{prop-qBexpdecay}, Proposition~\ref{prop-Cunionbound} can be concluded as follows.

\smallskip \noindent \emph{Proof of Proposition~\ref{prop-Cunionbound}.}
Let $L \in \N$ and let $z_1,\ldots,z_L \in \Z^d$ be pairwise distinct. By the stationarity of $\Lambda$, $ C_{n,M}(z_i)$, $i=1,\ldots,L$, are identically distributed. Using Lemma~\ref{lemma-mygoodness} and Proposition~\ref{prop-qBexpdecay}, we obtain for sufficiently large $n$, $\lambda=\lambda(n)$ and $M=M(\lambda,n)$ that
\begin{align*}
 \mathbb P & (C_{n,M}(z_1)^{\rm c} \cap \ldots \cap C_{n,M}(z_L)^{\rm c}) \\
& = \mathbb P((A_n(z_1)\cap B_{n,M}(z_1))^{\rm c} \cap \ldots \cap (A_{n}(z_L)\cap B_{n,M}(z_L))^{\rm c}) \\ 
&  \leq \mathbb P\Big( \big( \bigcup_{S \subseteq [L] \colon |S| \geq L/2} \bigcap_{i \in S} A_n(z_i)^{\rm c} \big) \cup \big( \bigcup_{ S \subseteq [L] \colon |S| \geq L/2} \bigcap_{i \in S} B_{n,M}(z_i)^{\rm c} \Big) \\ 
&  \leq 2 \max \Big\lbrace \P \Big(  \bigcup_{S \subseteq [L] \colon |S| \geq L/2} \bigcap_{i \in S} A_n(z_i)^{\rm c}   \Big) + \P \Big(\bigcup_{ S \subseteq [L] \colon |S| \geq L/2} \bigcap_{i \in S} B_{n,M}(z_i)^{\rm c}  \Big) \Big\rbrace \\
&   \leq 2 \binom{L}{\lfloor L/2 \rfloor} \max \lbrace q_A^{L/2}, q_B^{L/2} \rbrace \leq 2 \times 2^L \max \lbrace \sqrt{q_A}^L, \sqrt{q_B}^L \rbrace. 
\end{align*}
Putting  $q_C=4 \max \lbrace \sqrt{q_A}, \sqrt{q_B} \rbrace$  and choosing $n,\lambda,M$ large enough, the proposition follows. 
\ProofEnde 

\begin{step}\label{step-SIRgraphalsopercolates}
Percolation in the SINR graph.\index{cluster!infinite!in the SINR graph}\index{percolation!SINR!for Cox point processes}\index{percolation process!discrete!relation to percolation in the SINR graph}
\end{step}

Now, let $n,\lambda,M$ be such that the process of $n$-good sites $z \in \Z^d$ such that $I_{6n}(nz) \leq M$ percolates. If $z$ is such a site, then $I(x) \leq M$ for all $x \in Q_{6n}(n z)$. Now, as in \cite[Section~3.3]{DF06} in the case of a different discrete model, for an $n$-good site $z$ such that $I_{6n}(nz) \leq M$ and for $X_i,X_j \in X^\lambda \cap Q_{6n}(nz)$ with $|X_i-X_j| \leq r$, we have
\[ \frac{\ell(|X_i-X_j|)}{N_0+\gamma \sum_{k \neq i,j} \ell(|X_k-X_j|)} \geq \frac{ \ell(r)}{N_0+\gamma  M}.\]
Choosing
\[ \gamma' = \frac{N_0}{M} \Big( \frac{\ell(r)}{\ell(r_{\mathrm B})} -1 \Big)>0, \numberthis\label{goodgammadef} \]
(where the inequality holds because $\upsilon_0<r<r_{\mathrm B}$), we have \[ \frac{ \ell(r)}{N_0+\gamma'  M} = \frac{ \ell(r_{\mathrm B})}{N_0}=\tau. \numberthis\label{estimatedgammalambda} \] 
Thus, for $\gamma \in (0,\gamma')$, any two Cox points of distance less than $r$ both lying within $Q_{6n}(nz)$ for an $n$-good site $z$ such that $I_{6n}(nz) \leq M$ are connected in $g_{(\gamma)}(X^\lambda)$. 

    Finally, similarly to \cite[Section~5.2]{HJC17}, we have the following. If there exists an infinite connected component $\mathcal C$ of $n$-good sites $z$ with $I_{6n}(nz) \leq M$, let $z,z' \in \mathcal C$ with $|z-z'|=1$. Then by property \eqref{second-ngood1} in the definition of $n$-goodness, there exist $X_i \in \Xi_n(nz), X'_i \in \Xi_n(nz')$. By property \eqref{third-ngood1}, we find a path from $X_i$ to $X'_i$  in $g_r(X^\lambda) \cap Q_{6n}(nz)$. Since $I_{6n}(nz) \leq M$, all the edges of this path also exist in $g_{(\gamma)}(X^\lambda)$. Hence, $g_{(\gamma)}(X^\lambda) \cap (\bigcup_{z \in \mathcal C} Q_{6n}(nz))$ contains an infinite path, which implies that $g_{(\gamma)}(X^\lambda)$ percolates.

Thus, Theorem~\ref{prop-firstpercolation} follows as soon as we have proven Proposition~\ref{prop-qBexpdecay}. In Section~\ref{sec-additionalchallenges} we will comment on the arguments of this proof and possible generalizations in least technical terms.
\begin{step}\label{step-interferencecontrol}\index{interference!control|textbf}
Proof of Proposition~\ref{prop-qBexpdecay}.
\end{step}
We start the proof with splitting the interference into two parts. For $x \in \R^d$ and $n \geq 1$, we put
\begin{align*}
& I_{6n}^{\mathrm{in}} (x) = \sum_{X_i \in X^\lambda \cap Q_{12n\sqrt d}(x)} \ell_{6n}(|X_i-x|), \quad I_{6n}^{\mathrm{out}} (x) = \sum_{X_i \in X^\lambda \setminus Q_{12n\sqrt d}(x)} \ell_{6n}(|X_i-x|).
\end{align*}
Further, for $z \in \Z^d$, we write $B_{n,M}^{\mathrm{in}}(z)=  \lbrace I_{6n}^{\mathrm{in}}(nz) \leq M \rbrace$ and $B_{n,M}^{\mathrm{out}}(z)=\lbrace I_{6n}^{\mathrm{out}}(nz)$ $ \leq M \rbrace$. In order to verify Proposition~\ref{prop-qBexpdecay}, we will first prove the following assertion about the inner part of the interference.
\begin{proposition}\label{prop-qBinexpdecay}
Under the assumption \eqref{ellboundedcond} or \eqref{expmomentscond} in Theorem~\ref{prop-firstpercolation}, for all sufficiently large $n \geq 1$, for all $\lambda>0$, and for all sufficiently large $M=M(n,\lambda)>0$, there exists a constant $q_B=q_B(n,\lambda,M) <1$ such that for any $L \in \N$ and pairwise distinct sites $z_1,\ldots,z_L\in \Z^d$,
\[ \mathbb P(B^{\mathrm{in}}_{n,M}(z_1)^{\rm c} \cap \ldots \cap B_{n,M}^{\mathrm{in}}(z_L)^{\rm c}) \leq q_B^L. \numberthis\label{qBinexpdecay} \]
Moreover, for all $\eps>0$, for all large enough $n \geq 1$, and for all $\lambda>0$, we can choose $M$ large enough  such  that $q_B \leq \eps$.
\end{proposition}
Then we will show that the following proposition holds about the outer part of the interference under the exponential moment assumption in \eqref{expmomentscond} on $\Lambda$ for unboundedly supported $\ell$ satisfying Assumption~$(\ell)$.
\begin{proposition}\label{prop-qBoutexpdecay}
Under the assumption \eqref{expmomentscond} in Theorem~\ref{prop-firstpercolation} for unboundedly supported $\ell$, for all sufficiently large $n \geq 1$, for all $\lambda>0$, and for all sufficiently large $M=M(n,\lambda)>0$, there exists a constant $q_B=q_B(n,\lambda,M) <1$ such that for any $L \in \N$ and pairwise distinct sites $z_1,\ldots,z_L\in \Z^d$,
\[ \mathbb P(B^{\mathrm{out}}_{n,M}(z_1)^{\rm c} \cap \ldots \cap B_{n,M}^{\mathrm{out}}(z_L)^{\rm c}) \leq q_B^L. \numberthis\label{qBoutexpdecay} \]
Moreover, for all $\eps>0$, for all large enough $n \geq 1$, and for all $\lambda>0$, we can choose $M$ large enough  such  that $q_B \leq \eps$.
\end{proposition}
Before verifying these propositions, let us show how they imply Proposition~\ref{prop-qBexpdecay}. 
\begin{proof}[Proof of Proposition~\ref{prop-qBexpdecay}.]
Note that for $n,\lambda,M>0$, if $I_{6n}(x)>M$, then $I_{6n}^{\mathrm{in}} (x)>M/2$ or $I_{6n}^{\mathrm{out}} (x)>M/2$. Using a union bound, it suffices to verify Proposition~\ref{prop-qBexpdecay} both with $B_{n,M}(z_i)$ replaced by $B_{n,M/2}^{\mathrm{in}}(z_i)$ and with $B_{n,M}(z_i)$ replaced by $B_{n,M/2}^{\mathrm{out}}(z_i)$ everywhere in \eqref{qBexpdecay} for all $i \leq L$. Indeed, having these assertions, we can combine them analogously to the proof of Proposition~\ref{prop-Cunionbound}. Clearly, it is enough to prove them without the factors $1/2$ in front of $M$, since $M$ can be chosen arbitrary large in Proposition~\ref{prop-qBexpdecay}.

We conclude that under the assumption \eqref{expmomentscond} for unboundedly supported $\ell$, Propositions~\ref{prop-qBinexpdecay} and \ref{prop-qBoutexpdecay} imply Proposition~\ref{prop-qBexpdecay}. It remains to show that under the assumption \eqref{ellboundedcond} that $\ell$ has bounded support, Proposition~\ref{prop-qBinexpdecay} alone implies Proposition~\ref{prop-qBexpdecay}. But this is true because for $\ell$ compactly supported, for all sufficiently large $n$, the following holds for all $z \in \Z^d$
\[
\begin{aligned}
I_{6n}^{\mathrm{in}}(nz) \leq I_{6n}(nz) &= \sum_{X_i \in X^\lambda \cap Q_{6n + 2 \sup \supp(\ell)}(nz)} \ell_{6n}(|X_i-nz|) 
\\ & \leq \sum_{X_i \in X^\lambda \cap Q_{12n\sqrt d}(nz)} \ell_{6n}(|X_i-nz|) = I^{\mathrm{in}}_{6n}(nz).
\end{aligned}\]
\end{proof}
We now prove Proposition~\ref{prop-qBinexpdecay}.
\begin{proof}[Proof of Proposition~\ref{prop-qBinexpdecay}]
In order to carry out the proof, we fix $\lambda>0$ and construct a renormalized percolation process as follows. A site $z \in \Z^d$ is \emph{$n$-tame} if
\begin{enumerate}
\item\label{first-tame} $R(Q_{12n\sqrt d}(nz))<n/2$,
\item\label{second-tame} $I_{6n}^{\mathrm{in}}(nz)\leq M $.
\end{enumerate}
A site $z \in \Z^d$ is \emph{$n$-wild} if it is not $n$-tame. The process of $n$-tame sites is $\lceil 12 n \sqrt d + 1\rceil$-dependent according to the definition of stabilization. Thus, using\index{dependent percolation theory} dependent percolation theory \cite[Theorem~0.0]{LSS97} (similarly to the proof of Lemma~\ref{lemma-mygoodness}), in order to verify Proposition~\ref{prop-qBinexpdecay}, it suffices to show that $\mathbb P(z\text{ is }n\text{-wild})$ can be made arbitrarily close to 0 uniformly in $z \in \Z^d$ by choosing first $n$ sufficiently large in order that the condition \eqref{first-tame} is satisfied (note that this condition does not depend on $\lambda$), and then choosing $M=M(n,\lambda)$ large enough so that \eqref{second-tame} also holds. We have 
\[ \mathbb P(z\text{ is }n\text{-wild}) \leq \mathbb P(R(Q_{12n\sqrt d}(nz)) \geq n/2) + \mathbb P(I_{6n}^{\mathrm{in}}(nz)> M). \]
The first term can be made arbitrarily small by choosing $n$ large enough, according to the definition of stabilization. Further, by \eqref{elladef},
\[ I_{6n}^{\mathrm{in}}(nz) = \sum_{X_i \in X^\lambda \cap Q_{12n\sqrt d}(nz)} \ell_{6n}(|X_i-nz|) \leq \ell(0) \# \big( X^\lambda \cap Q_{12n\sqrt d}(nz) \big) \]
holds for all $z \in \Z^d$. In particular, 
\[ \mathbb E[I^{\mathrm{in}}_{6n}(nz)] \leq \ell(0) \lambda \E[\Lambda(Q_{12n\sqrt d})]=(12n\sqrt d)^d \ell(0) \lambda < \infty. \]
Thus, for any sufficiently large $n \geq 1$, $\mathbb P(I_{6n}^{\mathrm{in}}(nz)>M)$ can be made arbitrarily small uniformly in $z \in \Z^d$ by choosing $M=M(n,\lambda)$ large enough. Thus, we conclude Proposition~\ref{prop-qBinexpdecay}.
\end{proof}
It remains to prove Proposition~\ref{prop-qBoutexpdecay}. We start with a deterministic result about the shifted versions of path-loss functions $\ell$ satisfying Assumption~$(\ell)$, the use of which will become transparent during the proof of the proposition.
\begin{lemma}\label{lemma-K0}
If $\ell$ satisfies Assumption~$(\ell)$, then there exists $K_0>0$ such that for all $x \in \R^d$, $n \geq 1$, $L \in \N$ and pairwise distinct $z_1,\ldots,z_L \in \Z^d$ we have
\[ \sum_{i=1}^L \ell_{6n}(|x-nz_i|) \leq K_0. \]
\end{lemma}
\begin{proof} This proof follows \cite[Section 3.2]{DF06}, where only the case $d=2$ was considered. 
Under the assumptions of the lemma, since the sites $ z_i$, $i=1,\ldots,L$, are pairwise distinct, the sum $\sum_{i=1}^L \ell_{6n}(|x-n z_i|)$ can be upper bounded by $ \sum_{z \in \Z^d} \ell_{6n}(|x-n z|)$.
Further, the sites $ z_i$, $i=1,\ldots,L$, are contained in the hypercubic lattice $\Z^d$. Let us write $Q_x$ for the cube of $\Z^d$ containing $x$; this is well-defined for a.e.~$x \in \R^d$. Now, for such $x$, for $i \in \N_0$ and $z_0 \in \lbrace z \in \Z^d \colon i \leq \dist_\infty (z,Q_x)<(i+1) \rbrace$, the contribution of $\ell(|nz_0-x|)$ to the latter sum is at most $\ell(in)$. Thus, we have 
\[ 
\begin{aligned}
\sum_{i=1}^L \ell_{6n}(|x-n z_i|)  &\leq \sum_{z \in \Z^d} \ell_{6n}(|x-n z|) 
\\ &\leq \sum_{i=1}^{\infty} \# \lbrace z \in \Z^d \colon i \leq \dist_\infty (z,Q_x)<(i+1) \rbrace \ell_{6n}(i n)=: K(n), 
\end{aligned}\numberthis\label{interferenceupper} \]
where, \emph{a priori}, $K(n) \in [0,\infty]$. Since $\ell$ is decreasing, we have for any $n \geq 1$ that
\[ K(n) \leq 2^d + \sum_{i=0}^{\lceil 6\sqrt d/2\rceil} \ell(0) ((2i+2)^d-(2i)^d) + \sum_{i=\lceil 6\sqrt d/2\rceil}^{\infty} ((2i+2)^d-(2i)^d) \ell(i-6\sqrt d/2). \numberthis\label{K0def} \]
Writing $K_0$ for the expression on the right-hand side, by \eqref{last-pathloss}\index{path-loss function!integrability} in Assumption~($\ell$), we have that $K_0<\infty$, which implies the lemma.
\end{proof}
We now carry out the proof of Proposition~\ref{prop-qBoutexpdecay}. 
\begin{proof}[Proof of Proposition~\ref{prop-qBoutexpdecay}.]
We proceed similarly to \cite[Section 3.2]{DF06} (until the estimate \eqref{newstuff}). We fix $\lambda>0$. By\index{Markov's inequality!exponential} Markov's inequality, for any $s>0$,
\begin{align*}  \mathbb P(& B_{n,M}^{\mathrm{out}}(z_1)^{\rm c} \cap \ldots \cap B_{n,M}^{\mathrm{out}}(z_L)^{\rm c} )  = \mathbb P( I_{6n}^{\mathrm{out}}(nz_1)>M,\ldots,I_{6n}^{\mathrm{out}}(nz_L)>M ) 
\\  & \leq  \mathbb P \Big( \sum_{i=1}^L I_{6n}^{\mathrm{out}}(nz_i) > L M \Big) 
\\ & \leq  \e^{-sLM} \E \Big[ \exp\Big( s \sum_{i=1}^L \sum_{X_k \in X^\lambda \setminus Q_{12n\sqrt d}(nz_i)} \ell_{6n}(|nz_i-X_k|) \Big) \Big].\numberthis\label{Csebisev} \end{align*}
Applying the form of the\index{Laplace functional of Cox point processes}\index{point process!Cox}\index{Campbell's theorem!for Cox point processes} Laplace functional of a Cox point process (cf.~\cite[Sections~3.2,~6]{K93}) to the function $f(x)=s\sum_{i=1}^L \ell_{6n}(|x-z_i|) \ind{1}\lbrace x \in \R^d \setminus Q_{12n\sqrt d}(nz_i) \rbrace$, we obtain 
\begin{align*} &\E \Big[ \exp\Big( s \sum_{i=1}^L \sum_{X_k \in X^\lambda \setminus Q_{12n\sqrt d}(nz_i)}  \ell_{6n}(|nz_i-X_k|) \Big) \Big] \\ = &\E \Big[ \exp \Big( \lambda   \int_{\R^d } \Big[ \exp \Big( s \sum_{i=1}^L  \ell_{6n}(|nz_i-x|) \ind{1}\big\lbrace x \in \R^d \setminus Q_{12n\sqrt d}(n z_i) \big\rbrace  \Big) -1 \Big) \Lambda(\d x) \Big] \Big]. \numberthis\label{CoxLaplace} \end{align*}
Now, we provide a uniform upper bound on the sum on the right-hand side of \eqref{CoxLaplace}.
We fix $K_0$ satisfying the assumption of Lemma~\ref{lemma-K0} for the rest of the proof. Thus, choosing $s \leq 1/K_0$ in \eqref{CoxLaplace}, we see that $ s \sum_{i=1}^L \ell_{6n}(|nz_i-x|) \ind{1}\lbrace x \in \R^d \setminus Q_{12n\sqrt d}(n z_i) \rbrace \leq 1$. Therefore, using that $\exp(y)-1 \leq 2y$ for all $y \leq 1$, we have
\[
\begin{aligned}\exp& \Big( s \sum_{i=1}^L \ell_{6n}(|nz_i-x|) \ind{1}\big\lbrace x \in \R^d \setminus Q_{12n\sqrt d}(n z_i) \big\rbrace \Big)-1 \\ \leq 2s &\sum_{i=1}^L \ell_{6n}(|n z_i-x|) \ind{1}\big\lbrace x \in \R^d \setminus Q_{12n\sqrt d}(n z_i) \big\rbrace.\end{aligned}\numberthis\label{doubleexptrick} \]
Plugging this into \eqref{CoxLaplace}, we obtain
\[ 
\begin{aligned}
\E & \Big[ \exp\Big( s \sum_{i=1}^L \sum_{X_k \in X^\lambda \setminus Q_{12n\sqrt d}(nz_i)}  \ell_{6n}(|nz_i-X_k|) \Big) \Big] \\ & \leq 
\E \Big[  \exp \Big( 2\lambda s \sum_{i=1}^L \int_{\R^d \setminus Q_{12n\sqrt d}(n z_i)} \ell_{6n}(|nz_i-x|)\Lambda(\d x) \Big)
\Big].
\end{aligned}
\numberthis\label{newstuff} \]
We now estimate the right-hand side of \eqref{newstuff}. For $i \in [L]$, we extend the integration domains $\R^d \setminus Q_{12n\sqrt d}(nz_i)$ to $\R^d \setminus Q_{\lfloor 12n\sqrt d \rfloor}(nz_i)$ and we  
subdivide $\R^d \setminus Q_{\lfloor 12n\sqrt d \rfloor}(nz_i)$ into a union of concentric $\ell^\infty$-annuli $Q_{\lfloor 12n\sqrt d \rfloor+2}(nz_i) \setminus Q_{\lfloor 12n\sqrt d \rfloor}(nz_i)$, $Q_{\lfloor 12n\sqrt d \rfloor+4}(nz_i)\setminus Q_{\lfloor 12n\sqrt d \rfloor+2}(nz_i)$  etc.~(up to the boundaries). Now for each $j \in \N_0$, let us write $A^j= Q_{\lfloor 12n\sqrt d \rfloor +2j+2} \setminus Q_{\lfloor 12n\sqrt d \rfloor+2j}$. Note that $A^j$ is covered by the union of $\nu_j=(\lfloor 12n\sqrt d \rfloor+2j+2)^d  - (\lfloor 12n\sqrt d \rfloor +2j)^d $ congruent copies $Q^{j,1},Q^{j,2},\ldots,Q^{j,\nu_j}$ of $Q_1$, and $\nu_j=|A^j|$. Further, for $i \in [L]$ and $x \in A^j+nz_i$, we have for all sufficiently large $n$ (not depending on $j$)
\[ \ell_{6n}(|nz_i-x|) \leq \ell_{6n} \Big( \frac{\lfloor 12 n\sqrt d \rfloor+2j}{2} \Big) \leq \ell \Big( j-1+\frac{6n\sqrt d}{2} \Big) \leq \ell ( j+2n\sqrt d ) . \]
Hence, we obtain
\begin{align*}
 \mathbb E& \Big[  \exp \Big( 2 \lambda s \sum_{i=1}^L \int_{\R^d \setminus Q_{12n\sqrt d}(nz_i)}  \ell_{6n}(|nz_i-x|)\Lambda(\d x) \Big) \Big] 
\\ & \leq \mathbb E \Big[  \exp \Big( 2 \lambda s \sum_{i=1}^L \sum_{j=0}^{\infty} \Lambda \big( A^j+nz_i \big)  \ell (j+2n\sqrt d) \Big) \Big]
\\ & =  \mathbb E \Big[  \exp \Big( 2 \lambda s  \sum_{j=0}^{\infty} \sum_{k=1}^{\nu_j} \sum_{i=1}^L \Lambda \big( Q^{j,k}+nz_i \big)  \ell (j+2n\sqrt d) \Big) \Big].
\end{align*}
Since $\Lambda$ is $b$-dependent, for any $i,i' \in [L]$, $j \in \N_0$, and $k \in [\nu_j]$, the identically distributed random variables $\Lambda(Q^{j,k}+nz_i)$ and $\Lambda(Q^{j,k}+nz_{i'})$ are independent if $|nz_i-nz_{i'}| > b+\sqrt d$. In particular, since $n \geq 1$, $\Lambda(Q^{j,k}+nz_i)$ and $\Lambda(Q^{j,k}+nz_{i'})$ are $b+\sqrt d$-dependent. Hence, there exists a subset of $[L]$ with cardinality at least $\smfrac{L}{(b+\sqrt d+2)^d}$ such that $\Lambda(Q^{j,k}+nz_i)$ and $\Lambda(Q^{j,k}+nz_{i'})$ are independent for all $n \geq 1$, $j \in \N_0$ and $k \in[\nu_j]$. Hence, an application of Hölder's inequality yields
\begin{align*}
\mathbb E & \Big[  \exp \Big( 2 \lambda s  \sum_{j=0}^{\infty} \sum_{k=1}^{\nu_j} \sum_{i=1}^L \Lambda \big( Q^{j,k}+nz_i \big)  \ell (j+2n\sqrt d) \Big) \Big] \\
 \leq &  \mathbb E \Big[  \exp \Big(2 \lambda s (b+\sqrt d+2)^d \sum_{j=0}^{\infty} \sum_{k=1}^{\nu_j} \Lambda \big( Q^{j,k} \big)  \ell (j+2n\sqrt d) \Big) \Big]^{\frac{L}{(b+\sqrt d+2)^d}}. \numberthis\label{afterbdep}
\end{align*}

Using that the following extended version of Hölder's inequality holds for any sequence $(Y_i)_{i=1}^{\infty}$ of identically distributed non-negative random variables
\[ \mathbb E \Big(\prod_{i=1}^\infty  Y_i^{p_i} \Big) \leq \mathbb E (Y_1), \qquad  p_i \geq 0, ~\forall i \in \N,~\sum_{i=1}^{\infty} p_i=1, \numberthis\label{iteratedHoelder} \]
we obtain
\begin{align*}
\mathbb E & \Big[  \exp \Big(2 \lambda s (b+\sqrt d+2)^d \sum_{j=0}^{\infty} \sum_{k=1}^{\nu_j} \Lambda \big( Q^{j,k} \big)  \ell (j+2n\sqrt d) \Big) \Big] \\
\leq & \Big[  \exp \Big(2 \lambda s (b+\sqrt d+2)^d \sum_{j=0}^{\infty} |A_j|  \ell (j+2n\sqrt d) \Lambda(Q_1) \Big) \Big] . \numberthis\label{HoelderIntermediate}
\end{align*} 
Now, since for $i \in \N_0$, $|A^i| \leq  2  d (\lfloor 12 n \sqrt d \rfloor + 2i+2)^{d-1}$, for all sufficiently large $n$, the right-hand side of \eqref{HoelderIntermediate} is upper bounded by
\begin{align*}  \mathbb E &\Big[  \exp \Big(  4  \lambda s d \sum_{j=0}^{\infty} \big( 2+2j+ 7\times 2n\sqrt d \big)^{d-1} \ell \Big( j+2n \sqrt d \Big)  \Lambda(Q_1) \Big)\Big] \\ 
\leq \mathbb E & \Big[ \exp \Big(   4  \times 7^{d-1} \lambda s  d  \sum_{j=0}^{\infty} \big(j + 2n\sqrt d \big)^{d-1} \ell \Big( j+2n \sqrt d  \Big) \Lambda(Q_1) \Big)  \Big].
\end{align*}
Note that here, the requirement about how large $n$ has to be chosen does not depend on $\lambda$.

Now, $c_o=\sup_{n \geq 1} \sum_{j=0}^{\infty} (j + 2n\sqrt d )^{d-1} \ell ( j+2n \sqrt d )$ is finite thanks to \eqref{third-pathloss} in Assumption~$(\ell)$. Let us now fix $\alpha_o>0$ such that $\E[\exp(\alpha\Lambda(Q_1))]<\infty$, such $\alpha_o$ exists due to the exponential moment condition in \eqref{expmomentscond}. Now, it follows that given that $0<s \leq \max \{ \smfrac{1}{K_0}, \smfrac{\alpha_o}{4\times 7^{d-1} \lambda c_o d} \}$, we have
\[ 
\begin{aligned}
\mathbb P(& B_{n,M}^{\mathrm{out}}(z_1)^{\rm c} \cap \ldots \cap B_{n,M}^{\mathrm{out}}(z_L)^{\rm c} )
\\ \leq & \exp(-sML) \E[\exp(\alpha_o \Lambda(Q_1))]^{\frac{L}{(b+\sqrt d+2)^d}} \leq q_B^L,
\end{aligned} \numberthis\label{bdepdone}
\]
where $q_B:=\exp(-sM)\E[\exp(\alpha_o\Lambda(Q_1))]^{\frac{1}{(b+\sqrt d+2)^d}}$ can be made arbitrarily close to zero by choosing $M=M(n,\lambda)$ sufficiently large. This finishes the proof of Proposition~\ref{prop-qBexpdecay}.
\end{proof} 
\subsubsection{Discussion about the interference control argument}\label{sec-additionalchallenges}
Let us provide a number of remarks about the essential points of the proof of Theorem~\ref{prop-firstpercolation} part \eqref{firstsupercritical} and about its possible generalizations.
\begin{itemize}
    \item The proof of Proposition~\ref{prop-qBexpdecay} requires no asymptotic essential connectedness of $\Lambda$, only stabilization (and $b$-dependence under the condition \eqref{expmomentscond}). Asymptotic essential connectedness is only needed in order to be able to use strong connectivity of the underlying Gilbert graph, i.e., to verify Lemma~\ref{lemma-mygoodness} for $r_{\rm B}$ (and thus $r$) arbitrarily small (depending on the value of $d_0$).
    \item The proof of Lemma~\ref{lemma-K0} together with the integrability of $\ell$ shows that under the exponential-moment and $b$-dependence assumption \eqref{expmomentscond}, Proposition~\ref{prop-qBoutexpdecay} with $I_{6n}^{\text{out}}$ replaced by $I_{6n}$ everywhere can be proven along the lines of the proof of the original proposition. The reason for splitting the interference is just to be able to cover the case of only asymptotically essentially connected $\Lambda$ for $\ell$ compactly supported, since for such $\Lambda$ the proof of Proposition~\ref{prop-qBoutexpdecay} is not applicable.
    \item Roughly speaking, the final estimate \eqref{bdepdone} in the proof of Proposition~\ref{prop-qBoutexpdecay} tells that $B_{n,M}^{\mathrm{out}}(z_i)^{\rm c}$, $i \in [L]$, are $b'$-dependent for some $b'$ not depending on $L$. This holds despite the infinite-range spatial dependency of the interference measured at a given site $nz_i$, since exponential moments of this interference can be bounded by exponential moments of $\Lambda(Q_1)$ using Hölder's inequality and the integrability of $\ell$. Such an assertion seems difficult to show if $\Lambda$ is only stabilizing. Indeed, then even random variables of the form $\E[\exp(\alpha\Lambda(Q_1(nz_i)))]$, $i \in [L]$, are not $b$-dependent, and in the event that the suprema of stabilization radii are larger than $n$ in each of these $L$ boxes, exponential moments of $\sum_{i=1}^L \Lambda(Q_1(nz_i))$ can only be estimated by the ones of $L \Lambda(Q_1)$ (using Hölder's inequality). This grows superexponentially in $L$ and thus cannot be compensated by $\exp(-sML)$, unless $\Lambda(Q_1)$ is bounded. 
    Since we are not aware of any relevant intensity measure where $\Lambda(Q_1)$ is bounded and $\Lambda$ is stabilizing but not $b$-dependent, we omitted this case from the statement of Theorem~\ref{prop-firstpercolation} part \eqref{firstsupercritical}.
\end{itemize}

\subsection{Proof of the results of Section~\ref{sec-stabilizing}}\label{sec-stabilizingproof}\index{phase!supercritical!existence for stabilizing $\Lambda$ for large $r$}
\subsubsection{Proof of Corollary~\ref{cor-largersupercritical}}\label{sec-largersupercritical}
Before carrying out the proof, we recall\index{Palm!calculus} Palm calculus for\index{point process!Cox} Cox processes from \cite[Section~2.2]{HJC17}. The\index{Palm!version!of a stationary point process}\index{point process!stationary} \emph{Palm version} $X^{\lambda,*}$ of a stationary point process $X^\lambda$ of intensity $\lambda=\E[\# (X^\lambda \cap Q_1)]>0$ is a point process whose distribution is defined via
\[ \mathbb E\big[ f(X^{\lambda,*}) \big] = \frac{1}{\lambda} \mathbb E\Big[ \sum_{X_i \in X^\lambda \cap Q_1} f(X^\lambda-X_i) \Big],  \numberthis\label{Palmpointprocess} \]
for any bounded measurable function $f \colon \mathbb M_{\mathrm{co}} \to [0,\infty)$, where $\mathbb M_{\mathrm{co}}$ is the set of\index{measure!counting!$\sigma$-finite} $\sigma$-finite counting measures. In particular, $\mathbb P(o \in X^{\lambda,*})=1$. 

For any infinite, locally finite graph $G=(V,E)$ and for a vertex $v \in V$, we say that $v \leftrightsquigarrow \infty$ in $G$ if $v$ is contained in an infinite connected component of $G$. Then, for $r>0$,
\[ \theta(\lambda,r)=\mathbb P \Big( o \leftrightsquigarrow \infty \text{ in }g_{r}(X^{\lambda,*}) \Big) \numberthis\label{thetaincriticality} \]
denotes the\index{percolation probability!in the Cox--Gilbert graph} \emph{percolation probability} of the origin of  the Cox--Gilbert graph $g_r(X^\lambda)$. Then $\lambda_{\mathrm c}(r)=\inf \lbrace \lambda>0 \colon \theta(\lambda,r)>0 \rbrace$, cf.~\cite[Section~2.2]{HJC17}.

\noindent \emph{Proof of Corollary~\ref{cor-largersupercritical}.}\index{stabilization} We first verify \eqref{first-largersupercritical}. Let $\overline{\theta}(\varrho)$ be the\index{percolation probability!in the Poisson--Gilbert graph} percolation probability of the Gilbert graph of a stationary Poisson point process with intensity $\varrho>0$ and connection radius 1. By \cite[Theorem~2.9]{HJC17}, for $\Lambda$ stabilizing,
\[ \lim_{\begin{smallmatrix}& r \uparrow \infty, \lambda \downarrow 0, & \lambda r^d=\varrho \end{smallmatrix}} \theta(\lambda,r)=\overline{\theta}(\varrho). \]
Let $\varrho>0$ satisfy $\overline{\theta}(\varrho)>0$. If $\lambda r^d = \varrho$ for $r$ large enough, then $\theta(\lambda,r)>0$, thus $\lambda_{\mathrm c}(r) \leq \smfrac{\varrho}{r^d} < \infty$.  This verifies \eqref{first-largersupercritical}. Since $\smfrac{\varrho}{r^d} \to 0$ as $r \to \infty$, it follows also that $\lim_{r \to \infty} \lambda_{\mathrm c}(r)=0$. But this is \eqref{second-largersupercritical}. \ProofEnde

\subsubsection{Proof of Proposition~\ref{prop-stabilizingsupercritical}}\label{sec-stabilizingsupercritical} 
Although it is possible to provide one proof for all dimensions $d \geq 2$, we find it instructive to start with the case $d=2$ and to verify the assertion using crossings of $3n \times n$ boxes in that case. Indeed, this discrete model lead to the assertion that $\lambda_{N_0,\tau}=\lambda_{\mathrm c}(r_\mathrm B)$ in the two-dimensional Poisson case, and thus using this model may be helpful for future investigations of the precise value of $\lambda_{N_0,\tau}$ in the stabilizing Cox case for $d=2$ and large $r_\mathrm B$. Afterwards, we will sketch the proof for $d\geq 3$. 

Let $d =2$.
Let us write $\mathcal B(\widehat{\lambda},\widehat{r})$ for the\index{Boolean model!Poisson} Poisson--Boolean model with intensity $\widehat{\lambda}>0$ and connection radius $\widehat{r}>0$. Further, for $r>0$, let $\varrho_{\mathrm c}(r)$ be such that $\mathcal B(\varrho_{\mathrm c}(r),r)$ is critical. Then, by scale invariance\index{scale invariance!for Poisson--Boolean models} of Poisson--Boolean models, we have $\varrho_{\mathrm c}(r)=r^{-d}\varrho_{\mathrm c}(1)$. We fix $\varrho>\varrho_{\mathrm c}(1)$, then there exists $\varrho'<\varrho$ such that $\mathcal B(\varrho',1)$ is still supercritical. 

For $r>\upsilon_0$, let us write $r_{\mathrm B}(r)=\smfrac{\varrho}{\varrho'} r$ and $\lambda(r)=\varrho' r^{-d}$. Then by Assumption ($\ell$) \eqref{first-pathloss}, \eqref{third-pathloss}, and the fact that $\ell$ has unbounded support, $\ell(r_{\mathrm B}(r))<\ell(r)$ holds for all $r>\upsilon_0$. Further, let $N_0(r)$ and $\tau(r)$ be such that $r_{\mathrm B}(r)=\ell^{-1}(\tau(r)N_0(r))$; such parameters exists since $\ell$ has unbounded support and satisfies \eqref{first-pathloss} and \eqref{third-pathloss} in Assumption ($\ell$). We map the\index{Boolean model!Cox} Cox--Boolean model $\mathcal C(\lambda(r),r)=X^{\lambda(r)} \oplus B_{r/2}$ to a discrete edge percolation model similarly to \cite[Section~3.1]{DF06},\index{interference!control} control the interferences and conclude that if $r$ is large enough, then the SINR graph $g_{(\gamma,N_0(r),\tau(r))}(X^{\lambda(r)})$ with SNR connection radius $r_{\mathrm B}(r)$ percolates for some $\gamma>0$ (with probability 1 thanks to stabilization).

For $n \geq 1$ and $r >\upsilon_0$, let us write $z_e=(x_e,y_e)$ for the centre of the edge $e$ in the nearest neighbour graph of $\Z^2$. Let us denote the set of such edges by $E(\Z^2)$. Note that each $z_e$ is an element of $\mathbb X=\lbrace (x/2,y) \colon x,y \in \Z\rbrace \cup \lbrace (x,y/2) \colon x,y \in \Z \rbrace$. Let us write $T_e(n,r)=[nrx_e-\smfrac{3}{4}nr, nrx_e+\smfrac{3}{4}nr] \times [nry_e-\smfrac{1}{4}nr,nry_e+\smfrac{1}{4}nr]$ if $e \in E(\Z^2)$ is a horizontal edge and $T_e(n,r)=[nrx_e-\smfrac{1}{4}nr, nrx_e+\smfrac{1}{4}nr] \times [nry_e-\smfrac{3}{4}nr,nry_e+\smfrac{3}{4}nr]$ if $e$ is a vertical edge. Note that $S_e$ is a rectangle with its edges parallel to $e$ having length $\frac{3}{2}nr$ and its edges perpendicular to $e$ having length $\frac{1}{2}nr$. In particular, $Q_{nr/2}(nrz_e) \subset T_e(n,r) \subset Q_{3nr/2}(nrz_e)$, and $T_e(n,r) \setminus Q_{nr/2}(nrz_e)$ is the disjoint union of two $\smfrac{nr}{2} \times \smfrac{nr}{2}$ squares, let us denote their closures by $S_e^1(n,r)$ respectively $S_e^2(n,r)$ (in an arbitrary but fixed order for each $e$).
For any edge $e$ in $E(\Z^2)$, we say that $e$ is \emph{$(n,r)$-good} if\index{crossing probabilities}
\begin{enumerate}
\item\label{first-3nn} $R(Q_{\frac{3}{2}nr}(nrz_e))<\frac{3}{2}nr$, and
\item\label{second-3nn} $\mathcal C(\lambda(r),r)$ crosses $T_e(n,r)$ in the hard direction and both $S_{e}^1(n,r)$ and $S_e^2(n,r)$ in the other direction.
\end{enumerate}
An edge $e$ is \emph{$(n,r)$-bad} if it is not $(n,r)$-good. The process of $(n,r)$-good edges is\index{percolation process!discrete!$b$-dependent} 4-dependent as can be seen from the definition of stabilization. We write $J_{n,r}(z_e)$ for the event in \eqref{second-3nn} and $F_{n}(z_e)$ for the event that $\mathcal B(\varrho',1)$ crosses $T_e(n,1)$ in the hard direction and both $S_{e}^1(n,1)$ and $S_e^2(n,1)$ in the other direction. (The precise definitions of these events are analogous to the one in Theorem~\ref{theorem-MR96}, therefore we leave them to the reader.) Note that by scale invariance of the Poisson--Boolean model, $F_n(z_e)$ has probability equal to the one of the event that $\mathcal B(\lambda(r),r)=\mathcal B(\frac{\varrho'}{r^2},r)$ crosses $T_{e'}(n,r)$ in the hard direction and both $S_{e'}^1(n,r)$ and $S_{e'}^2(n,r)$ in the other direction for an arbitrary $e' \in E(\Z^2)$ and $r>0$. 

Now, let $\eps>0$. First, we fix $n$ sufficiently large  such  that the probability of the event in \eqref{first-3nn} is at least $1-\eps/4$ uniformly for all $r \geq 1$ and any edge $e$ in $n r \mathbb Z^2$, and that the probability of $F_{n}(z_e)^{\mathrm c}$ is also at most $\eps/4$. The last condition can be satisfied thanks to the\index{Russo--Seymour--Welsh type result} Russo--Seymour--Welsh type result, Theorem~\ref{theorem-MR96}. Next, as observed in \cite[Section~7.1]{HJC17}, the restriction of $r^{-1} \mathcal C(\lambda(r),r)$ to a bounded sampling window converges weakly to the corresponding restriction of $\mathcal B(\varrho',1)$. Now, for fixed $e$, the event $F_n(z_e)$ has discontinuities of measure 0 with respect to the Poisson--Boolean model. This implies that for all $r>\upsilon_0$ sufficiently large, $|\mathbb P(F_{n}(z_e)^{\mathrm c})-\mathbb P(J_{n,r}(z_e)^{\mathrm c})|$ can be bounded from above by $\eps/4$ uniformly in $e$. Now, for any $e$, using a union bound and the triangle inequality, we have
\[ \begin{aligned}\mathbb P&(e\text{ is }(n,r)\text{-bad}) \\ &\leq \mathbb P \Big( R(Q_{\frac{3}{2}rn}(z_e)) \geq \frac{3}{2}nr \Big) +\mathbb P(F_n(z_e)^{\mathrm c})+|\mathbb P(F_{n}(z_e)^{\mathrm c})-\mathbb P(J_{n,r}(z_e)^{\mathrm c})| \leq \frac{3\eps}{4}<\eps.\end{aligned} \]
Applying\index{dependent percolation theory}\index{domination by product measures} \cite[Theorem~0.0]{LSS97}, for all sufficiently large $n$ and large enough $r$ chosen accordingly, the process of $(n,r)$-good edges is\index{domination by product measures} stochastically dominated from below by\index{percolation process!discrete!edge (bond)} a supercritical independent edge percolation process. Thus, the $(n,r)$-good sites percolate for all sufficiently large $n,r$.

Next, the\index{interference!control} interferences can be controlled analogously to Proposition~\ref{prop-qBexpdecay}. Instead of $\lbrace I_{6n}(nz)\colon z \in \Z^d\rbrace$ in Step~\ref{step-latticepercolation} defined in Section~\ref{sec-firstpercolationproof}, now one should work with the rescaled interferences $\lbrace I_{3rn/2}(nrz_e) \colon e \in E(\Z^2)\rbrace$ associated to the edges. For $n,r \geq 1$, $M>0$, and $e \in E(\Z^2)$, let us write $B_{n,r,M}(e)$ for the event that $I_{3rn/2}(nrz_e) \leq M$. Under the assumption \eqref{expmomentscond}, it can be proven analogously to Proposition~\ref{prop-qBexpdecay} that for any $L \in \N$ and for any pairwise distinct $e_1,\ldots,e_L$, 
\[ \mathbb P(B_{n,r,M}(e_1)^{\rm c} \cap \ldots \cap B_{n,r,M}(e_L)^{\rm c}) \leq q_B^L  \]
for some $q_B \in [0,1)$, where for fixed, large enough $n,r$ and $\lambda=\lambda(r)$, $q_B$ can be made arbitrarily close to 0 by choosing $M$ sufficiently large. Using a\index{Peierls argument} Peierls argument, we see that for all sufficiently large $n$, $r$ (depending on $n$), and $M$ (depending on $n,r$), the process of $(n,r)$-good edges $e$ with $I_{3rn/2}(nrz_e) \leq M$ percolates. Just as in \cite[Sections~3.2,~3.3]{DF06}, this\index{percolation process!discrete!relation to percolation in the SINR graph} implies percolation of $g_{(\gamma,N_0(r),\tau(r))}(X^{\lambda(r)})$ for $\gamma \in (0,\gamma^*(r))$, where 
\[ \gamma^*(r) = \frac{N_0(r)}{M} \Big( \frac{\ell(r)}{\ell(r_{\mathrm B}(r))} - 1 \Big) >0 \] 
(cf.~\eqref{estimatedgammalambda}, here we used again that $r_{\mathrm B}(r)>r>\upsilon_0$). This holds whenever $r_{\mathrm B}(r)=\ell^{-1}(\tau(r)N_0(r))$. Thus, since $\lambda(r) \downarrow 0$ as $r \to \infty$, Proposition~\ref{prop-stabilizingsupercritical} follows for small enough $\lambda>0$. But increasing $\lambda$ increases the probability of $(n,r)$-goodness of any edge, and it is easy to see that also the analogue of Proposition~\ref{prop-qBexpdecay} works for larger $\lambda>0$ (at the price of reducing $\gamma>0$ without vanishing). We conclude Proposition~\ref{prop-stabilizingsupercritical} for $d=2$.

For $d \geq 3$, one can proceed with an analogous definition of all the parameters from the first two paragraphs of the proof for $d=2$ (adapted to the value of $d$), using a different discrete model. Here, we shall define a site $z \in \Z^d$ to be $(n,r)$-good if it satisfies the definition of $n$-goodness in Section~\ref{sec-firstpercolationproof}, Step~\ref{step-latticemapping}, but with $n$ replaced by $nr$ and $\lambda$ by $\lambda(r)$ (in particular, with $X^{\lambda}$ replaced by $X^{\lambda(r)}$ and $g_r(X^\lambda)$ replaced by $g_r(X^{\lambda(r)})$) everywhere. For $z \in \Z^d$, we write $J_{n,r}(z)$ for the event that $z$ satisfies \eqref{second-ngood1} and \eqref{third-ngood1} in the definition of $(n,r)$-goodness. Then, for any $n,r$ under consideration, the process of $(n,r)$-good sites is 7-dependent according to the definition of stabilization. Further, we let $Y^{\varrho'}$ be a Poisson point process with intensity ${\varrho'}=\lambda(r) r^d$, and we write $F_n(z)$ for the event that in the definition of $(n,1)$-goodness, $z$ satisfies \eqref{second-ngood1} with $X^\lambda$ replaced by $Y^{\varrho'}$, and \eqref{third-ngood1} with $g_r(X^\lambda)$ replaced by $g_1(Y^{\varrho'})$ everywhere. The probability of $F_n(z)$ is independent of the choice of $z$ and tends to one as $n \to \infty$ thanks to the arguments of \cite[Section 5.2]{HJC17}, since the constant intensity measure of the Poisson point process $Y^{\varrho'}$ is obviously asymptotically essentially connected. Using  the \index{scale invariance!for Poisson--Gilbert graphs}scale invariance of Poisson--Gilbert graphs, we conclude that for $z \in \Z^d$,
\[ \mathbb P(z\text{ is }(n,r)\text{-bad}) \leq \mathbb P \Big( R(Q_{6nr}(nrz)) \geq \frac{nr}{2} \Big) +\mathbb P(F_n(z)^{\mathrm c})+|\mathbb P(F_{n}(z)^{\mathrm c})-\mathbb P(J_{n,r}(z)^{\mathrm c})|, \]
which can be made arbitrarily large by first choosing $n$ large and then $r$ large according to $n$, thanks to the weak convergence of $r^{-1} g_r(X^\lambda)$ to $g_1(Y^{\varrho'})$ as $r \to \infty$, $\lambda(r) \to 0$, $r^d \lambda(r) = {\varrho'}$. Thus, the proof for\index{percolation!SINR!in higher dimensions} $d \geq 3$ can be completed analogously to the case $d=2$, and the scale invariance of Poisson--Gilbert graphs also implies that $\lambda_{N_0(r),\tau(r)}$ tends to zero in this coupled limit. As already indicated in Section~\ref{sec-stabilizing}, the proof for $d\geq 3$ is also applicable for $d=2$. 
\ProofEnde

\section{Proof of results about the critical interference cancellation factor}\label{sec-gammalambdaproof}
In this section we prove the results of Section~\ref{sec-gammalambda}. Section~\ref{sec-lambdadependentgammaproof} contains the proofs of the results of Section~\ref{sec-lambdadependentgammabound}, in particular we prove Proposition~\ref{prop-gammalambda0} in Section~\ref{sec-gammalambda0proof} and Proposition~\ref{prop-bdependentOlambda} in Section~\ref{sec-bdependentOlambdaproof}. In Section~\ref{sec-Omega(1/lambda)proof} we sketch the proof of Corollary~\ref{cor-Omega1/lambda}, the result of Section~\ref{sec-Omega(1/lambda)discussion}.

Throughout this section, we will use the simplified notation $g_{(\gamma)}(X^\lambda)$ instead of the one $g_{(\gamma,N_0,\tau)}(X^\lambda)$, and similarly for directed SINR graphs, because the parameters $N_0$ and $\tau$ are fixed in the context of these proofs.
\subsection{Proof of the results of Section~\ref{sec-lambdadependentgammabound}}\label{sec-lambdadependentgammaproof}
\subsubsection{Proof of Proposition~\ref{prop-gammalambda0}}\label{sec-gammalambda0proof}\index{interference cancellation factor!estimates on the critical value!upper bound for large $\lambda$}
We first consider the case $N_0>0$. It suffices to show that for fixed $N_0,\tau,\gamma>0$, there is $\lambda_0>0$ such that for all $\lambda>\lambda_0$, $\P(g_{(\gamma)}(X^\lambda)\text{ percolates})=0$. By Proposition~\ref{prop-gammalambdaestimate1}, the statement is clear if $\gamma \geq \frac{1}{\tau}$. Else, let $L \geq 2$ be such that $\gamma \geq \frac{1}{(L-1)\tau}$. By \eqref{directeddegreebound}, all in-degrees in $g^{\rightarrow}_{(\gamma)}(X^\lambda)$ are at most $L-1$. Let now\index{percolation process!discrete!site} $(Q^{j})_{j=1}^{\infty}$ be a subdivision of $\R^d$ into congruent copies of $Q_{\upsilon_0'/\sqrt d}$, where $\upsilon_0'=\ell^{-1}(\ell(0)/2)$ exists by Assumption~($\ell$). Then, for any $j \in \N$, we have $\ell(|x-y|) \in [\ell(0)/2, \ell(0)]$ for all $x,y \in Q^j$. 

We claim that\index{percolation process!discrete!relation to percolation in the SINR graph} if $g_{(\gamma)}(X^\lambda)$ percolates, then each $Q^j$ containing at least one point $X_i\in X^\lambda$ from an unbounded cluster of  $g_{(\gamma)}(X^\lambda)$ contains at most $2L+2$ points of $X^\lambda$. Indeed, otherwise, since $X_i$ is not isolated in $g_{(\gamma)}(X^\lambda)$, there exists $k \neq i$ such that $X_k \to X_i$ is an edge in $g^{\rightarrow}_{(\gamma)}(X^\lambda)$.  Now, if at least $2L$ points of $X^\lambda \setminus \lbrace X_i,X_k \rbrace$ are within distance at most $\upsilon_0'$ from $X_i$, then 
\[ \SINR(X_k,X_i,X^\lambda)= \frac{\ell(|X_k-X_i|)}{N_0+\gamma \sum_{j \neq k,i} \ell(|X_j-X_i|)} \leq \frac{\ell(0)}{2L\gamma \frac{\ell(0)}{2}} \leq \tau, \]
where in the last step we have used that $\gamma L \geq \smfrac{L}{(L-1)\tau} > \frac{1}{\tau}$. This implies the claim.

Since $N_0>0$, any edge in $g_{(\gamma)}(X^\lambda)$  has length at most $r_{\mathrm B}$. Thus, if $g_{(\gamma)}(X^\lambda)$ percolates, then so does the process of open sites in the following site percolation model defined on the set of centres $C(Q^i)$ of the boxes $Q^i$, $i \in \N$. The site $C(Q^i)$, $i \in \N$, is \emph{open} if there exists $j \in \N$ such that $\# (Q^j \cap X^\lambda) \leq 2L+2$ and $\min_{x \in Q^i, y \in Q^j}|x-y| \leq r_{\mathrm B}$. Here, the edge set of the site percolation model corresponds to the $\ell^1$-neighbourhood of the sites.

We now show that the process of open sites does not percolate for $\lambda$ large, almost surely. This process is clearly $b'$-dependent for sufficiently large $b'>0$ because $X^\lambda$ is $b$-dependent, and openness of a site depends on points of $X^\lambda$ in a bounded neighbourhood of the site. Thus, it suffices to show that $\mathbb P(C(Q^i)\text{ is open})$ tends to zero as $\lambda \to \infty$ uniformly in $i$. Indeed, applying\index{dependent percolation theory} dependent percolation theory \cite[Theorem~0.0]{LSS97}, for large $\lambda$, the process of open sites is\index{domination by product measures} stochastically dominated by a\index{percolation process!independent Bernoulli!site} subcritical independent Bernoulli site percolation process. By stationarity of $X^\lambda$, for all $i$ we have the union bound
\begin{align*}
& \mathbb P( C(Q^i)\text{ is open}) \leq \mathbb P \Big( \exists j \colon \min_{x \in Q^j, y \in Q^i}|x-y| \leq r_{\mathrm B},~ \# (X^\lambda \cap Q^j) \leq 3L \Big) \\
\leq & C r_{\mathrm B}^d   \mathbb P(\# (X^\lambda \cap Q_{\upsilon_0'/\sqrt d}) \leq 3L) = C r_{\mathrm B}^d  \mathbb E \Big( \e^{-\lambda \Lambda(Q_{\upsilon_0'/\sqrt d})} \sum_{k=0}^{3N} \frac{\lambda \Lambda(Q_{\upsilon_0'/\sqrt d})^k}{k!} \Big),
\end{align*}
for a suitably large constant $C>0$. Clearly, the right-hand side tends to 0 as $\lambda \to \infty$.


The case that $\ell$ has bounded support (and possibly $N_0=0$) can be handled analogously, replacing $r_{\mathrm B}$ with $\sup\supp(\ell)$, which is a bound on the length of any edge in $g_{(\gamma)}(X^\lambda)$ in this case.
\ProofEnde

\subsubsection{Proof of Proposition~\ref{prop-bdependentOlambda}}\label{sec-bdependentOlambdaproof}\index{interference cancellation factor!estimates on the critical value!upper bound for large $\lambda$}\index{percolation!SINR!in two dimensions}  We fix $d=2$, $N_0,\tau>0$, and $M>\ell(0)$. Further, we fix $\delta>0$ and $c_0>0$ such that $\ell(r)>\tau N_0$ for all $r \in [0, \delta]$ and $\mathbb P(\Lambda(Q_{\delta/2})>c_0)=1$. 


The proof is based on \cite[Section~III-D]{DBT03} in the Poisson case. Let us summarize that proof in a way that is adaptable to the Cox case.  
The authors of \cite{DBT03} constructed\index{percolation process!discrete!site} a square lattice with edge length $\delta/2$ with $o$ being situated in the centre of a square. They showed that for any square of this lattice, if the number of Poisson points in the square is more than $L'=\smfrac{(1+2\tau\gamma)M}{\tau^2 \gamma N_0}>0$, then all Poisson points in this square are isolated in $g_{(\gamma)}(X^\lambda)$. This also holds if $X^\lambda$ is replaced by any stationary Cox process (or even by any \index{point process!simple} simple point process). Let us call a square \emph{open} if it has at most $2L'$ Poisson points and \emph{closed} otherwise.

Next, by the independence property of the\index{point process!Poisson} Poisson point process, any two squares are open or closed independently of each other, and thus the open sets form an\index{percolation process!independent Bernoulli!site} independent Bernoulli site percolation process. Now, by elementary properties of the\index{Poisson distribution} Poisson distribution \cite[Lemma 1]{DBT03}, this process is subcritical for all $\lambda$ sufficiently large, in which case the origin is almost surely surrounded by a\index{circuit of closed squares} circuit of closed squares. Then, the proof was concluded by verifying the following statement.\index{percolation process!discrete!relation to percolation in the SINR graph} If the origin is surrounded by a circuit of closed squares, then for any $X_i \in X^\lambda \cap Q_{\delta/2}$, we have $X_i \not \leftrightsquigarrow \infty$. Indeed, the statement is clear if $Q_{\delta/2}$ is itself a closed square. Else, as it was shown in \cite[Theorem~III-D]{DBT03}, if $X_i,X_j \in X^\lambda$ are situated on two different sides  of a circuit of closed squares, then $\SINR(X_i,X_j,X^\lambda) \leq \tau$.  This statement is entirely deterministic and remains true after replacing $X^\lambda$ with any stationary Cox (or even any simple stationary) point process. It follows that $\mathbb E[\# \lbrace X_i \in X^\lambda \cap Q_{\delta/2} \colon X_i \leftrightsquigarrow \infty\text{ in $g_{(\gamma)}(X^\lambda)$} \rbrace]=0$, and thus $\mathbb P(\text{$g_{(\gamma)}(X^\lambda)$ percolates})=0$ by stationarity. 

In the $b$-dependent Cox case, the process of closed sites is\index{percolation process!discrete!$b$-dependent} 
$b'$-dependent for all sufficiently large $b'$.
Our goal is to show that
\[ \lim_{\lambda \to \infty} \mathbb P(\text{a given square is closed}) = 1. \numberthis\label{Palmfett} \]
Having this, by \cite[Theorem~0.0]{LSS97}, the process of closed sites is stochastically dominated\index{domination by product measures} from below by a supercritical independent Bernoulli percolation process for large enough $\lambda$, and thus almost surely there exists an circuit of closed squares surrounding $o$. This allows us to conclude the proposition analogously to \cite[Section~III-D]{DBT03}.


Now we verify \eqref{Palmfett}. For $\mu>0$, we write $Y(\mu)$ for a Poisson random variable with mean $\mu$. Let $\eps \in (0,1)$. In order to simplify the notation we write $X^\lambda(\cdot) = \# (X^\lambda \cap \cdot)$. By\index{Chebyshev's inequality} Chebyshev's inequality, we have
\begin{align*}
& \mathbb P \Big( \big| X^\lambda(Q_{\delta/2}(x))- \lambda \Lambda(Q_{\delta/2}(x))  \big| >\eps \lambda \Lambda(Q_{\delta/2}(x)) \Big) 
\\ = & \mathbb E \Big(  \mathbb P \Big( \big|   Y(\lambda \Lambda(Q_{\delta/2}(x)))- \mathbb E [Y(\lambda \Lambda(Q_{\delta/2}(x))) \big| >\eps \mathbb \E [Y(\lambda \Lambda(Q_{\delta/2}(x)))] \Big| \Lambda \Big)  \Big)
\\ \leq  & \mathbb E \Big( \mathbb E \Big( \frac{\mathrm{Var}(Y(\lambda \Lambda(Q_{\delta/2})))}{\eps^2 \mathbb E [Y(\lambda \Lambda(Q_{\delta/2}(x)))]^2} \Big| \Lambda \Big) \Big) = \mathbb E \Big( \frac{\lambda \Lambda(Q_{\delta/2}(x))}{\eps^2 \lambda^2 \Lambda(Q_{\delta/2}(x))^2} \Big) 
\\  = &\frac{1}{\eps^2\lambda} \E \Big( \frac{1}{\Lambda(Q_{\delta/2}(x))}\Big).  \numberthis\label{towerPalm}
\end{align*}
Under the assumption that $\Lambda(Q_{\delta/2})>c_0$ almost surely, the right-hand side is finite for all $\lambda>0$ and tends to 0 as $\lambda \to \infty$. 
Now, similarly to \cite[Section~III-D]{DBT03}, if $\lambda$ satisfies
\[ \lambda \geq \frac{2L'}{(1-\eps) c_0}, \numberthis\label{howtochooselambda} \]
where, almost surely, the right-hand side is more than $\smfrac{L'}{(1-\eps)\Lambda(Q_{\delta/2}(x)) }$ for all $x \in \R^2$, then 
\begin{align*}
\mathbb P& \big( X^\lambda( Q_{\delta/2}(x)) \leq L' \big) \leq  \mathbb P \big( X^\lambda( Q_{\delta/2}(x))  \leq (1-\eps)c_0\lambda \big) 
 \\ \leq \mathbb P &\big( X^{\lambda}(Q_{\delta/2}) \leq (1-\eps)\lambda \Lambda(Q_{\delta/2}(x)) \big) 
\\ \leq   \mathbb P &\Big( \big|   X^{\lambda}(Q_{\delta/2}(x))- \lambda \Lambda(Q_{\delta/2}(x))  \big| >\eps \lambda \Lambda(Q_{\delta/2}(x)  \Big),
\end{align*}
and thus by \eqref{towerPalm}, \eqref{Palmfett} holds. By \eqref{howtochooselambda}, $\gamma>\gamma^*(\lambda)$ holds once
\[ \frac{2(1+2\tau\gamma)M}{\tau^2 \gamma N_0} \leq (1-\eps)\lambda c_0,\]
or equivalently,
\[ \gamma \geq \frac{2M}{(1-\eps)\tau^2 N_0 \lambda c_0-4\tau  M}. \]
This is true if
\[ \gamma
\geq \frac{M}{(1-2\eps)\tau^2 N_0 \lambda c_0}, \numberthis\label{howtochoosegamma} \]
for $\lambda$ sufficiently large, namely for $\lambda \geq \smfrac{4M}{\eps \tau N_0 c_0}$. Clearly, the lower bound on the right-hand side of \eqref{howtochoosegamma} is in $\Ocal(1/\lambda)$. Now, for $\lambda \geq \smfrac{4M}{\eps \tau N_0 c_0}$ so large that the process of closed sites is\index{domination by product measures} stochastically dominated from below by a supercritical independent Bernoulli percolation process, the origin is almost surely surrounded by a circuit of closed squares. We conclude the proposition. 
\ProofEnde

\subsection{Sketch of proof of Corollary~\ref{cor-Omega1/lambda}}\label{sec-Omega(1/lambda)proof}\index{interference cancellation factor!estimates on the critical value!lower bound for large $\lambda$}\index{percolation!SINR!in two dimensions}
Since the assertion of the corollary is a lower bound on $\gamma^*(\lambda)$ for large $\lambda$, it suffices to verify it for $N_0>0$ (cf.~Section~\ref{sec-noIse}). We fix $d=2$, $\tau,N_0>0$, $M>\ell(0)$, and $\delta>0$ such that $\ell(\delta)>\tau N_0$. Further, we assume that $\eta \in [\sup \supp(\ell),\infty)$ such that\index{measure!intensity!locally bounded away from zero} $\Lambda(Q_{\eta})$ is bounded away from 0 (such an $\eta$ exists by the assumption of Corollary~\ref{cor-Omega1/lambda}); we will make stronger assumptions on $\eta$ later during the proof.

In the following, we summarize the proof of the assertion $\gamma^*(\lambda)=\Omega(1/\lambda)$ from \cite[Section~III-C]{DBT03} in the Poisson case and afterwards we explain how it can be extended to the setting of Corollary~\ref{cor-Omega1/lambda}. For $\lambda>0$, one maps the SINR graph $g_{(\gamma)}(X^\lambda)$ to a square lattice $\mathcal H$ with edge length $\eta$. The\index{dual lattice} dual lattice of $\mathcal H$, i.e., $\mathcal H$ shifted by the vector $(\eta/2,\eta/2)$, is denoted by $\mathcal H'$. Lote that there is a one-to-one correspondence between the edges of $\mathcal H$ and the ones of $\mathcal H'$ by mapping an edge $e$ of $\mathcal H$ to the unique edge of $\mathcal H'$ which crosses $e$. In $\mathcal H$, one divides each square into $K^2$ subsquares of size $\smfrac{\eta}{K} \times \smfrac{\eta}{K}$, where $K \in \N$ is defined as $K=\lceil \smfrac{\sqrt 5 \eta}{\delta} \rceil$. Further, for $\lambda,\gamma>0$, one puts
\[ L =\inf_{x \colon |x| \leq \sqrt 5 \eta/K} \Big\lfloor \frac{1}{\gamma M} \Big( \frac{\ell(|x|)}{\tau}-N_0 \Big) \Big\rfloor=\Big\lfloor \frac{1}{\gamma M} \Big( \frac{\ell(\sqrt 5 \eta/K )}{\tau}-N_0 \Big) \Big\rfloor . \numberthis\label{Ndef} \]
One says that\index{percolation process!discrete!edge (bond)} a square of $\mathcal H$ is \emph{populated} if each of its subsquares contain at least one point of $X^\lambda$. Further, an edge $a$ of $\mathcal H$ is \emph{open} if both squares adjacent to $a$ are populated and the total number of points of $X^\lambda$ in all squares of $\mathcal H$ having at least one vertex in common with the neighbouring squares of $a$ is at most $L+1$. An edge $a'$ of $\mathcal H'$ is \emph{open} if and only if the corresponding edge of $\mathcal H$ is open. The proof proceeds by the following lemma \cite[Lemma 2]{DBT03}.
\begin{lemma}[\cite{DBT03}]\label{lemma-DBTLemma2}
Let $p$ denote the probability that an arbitrary edge in $\mathcal H'$ is closed, and let us write $q=1-p$. Then for any $q'>0$, there exists $\lambda' \in (0,\infty)$ such that for all $\lambda>\lambda'$ there exists $\gamma'(\lambda)>0$ such that for all $\gamma \in (0,\gamma'(\lambda)]$, $q<q'$. Further, $\lambda \mapsto \gamma'(\lambda)$ can be chosen such a way that $\gamma'(\lambda)=\Omega(1/\lambda)$ as $\lambda \to \infty$. 
\end{lemma}
The process of open edges in $\mathcal H'$ is 3-dependent thanks to the\index{point process!Poisson} independence property of the Poisson point process. Using\index{dependent percolation theory} dependent percolation theory \cite{LSS97}, one concludes that for all sufficiently large $\lambda>0$, there exists $\gamma'(\lambda)>0$ such that for all $\gamma \in (0,\gamma'(\lambda))$, the process of open edges percolates with probability 1, and such that $\gamma'(\lambda)=\Omega(1/\lambda)$. This\index{percolation process!discrete!relation to percolation in the SINR graph} implies percolation in the SINR graph thanks to \cite[Lemmas 4, 5]{DBT03}. These lemmas are similar to Step~\ref{step-SIRgraphalsopercolates} of the proof of Theorem~\ref{prop-firstpercolation}\eqref{firstsupercritical}, they use that $\ell$ has bounded support and $\eta \geq \sup \supp(\ell)$, but they are easily seen to hold for any\index{point process!simple} simple point process rather than only for the Poisson one.

Now, if $X^\lambda$ is a Cox point process with intensity $\lambda\Lambda$ where $\Lambda$ is $b$-dependent, then the process of open edges in $\mathcal H'$ is still $b'$-dependent for all sufficiently large $b'$. Thus, in order to conclude Corollary~\ref{cor-Omega1/lambda}, it suffices to verify Lemma \ref{lemma-DBTLemma2} under the assumptions of the corollary, for $\eta$ sufficiently large. 

\emph{Proof of Lemma~\ref{lemma-DBTLemma2} under the assumptions of Corollary~\ref{cor-Omega1/lambda}.} 
For $p=1-q$, we estimate
\begin{align*}
p  &=  \mathbb P(  \text{a given edge of $\mathcal H'$ is open}) \\
  & \geq   \mathbb P \big( \text{$2 K^2$ subsquares of area $\eta^2/K^2$ have at least 1 point of $X^\lambda$ each, and} 
  \\
  & \qquad \text{an area of $12 \eta^2$ including these subsquares has at most $N$ points} \big) \\
 & \geq   \mathbb P \big( \text{$12 K^2$ subsquares of area $\eta^2/K^2$ have between 1 and $\lfloor L/(12 K^2) \rfloor$} \\ & \qquad \text{points of $X^\lambda$ each} \big)
  \\ &= 1-\mathbb  P \big(  \text{at least 1 of $12 K^2$ subsquares of area $\eta^2/K^2$ has 0} \\& \qquad \text{or more than $\lfloor L/(12 K^2) \rfloor$ points}\big) 
  \\  &\geq 1- 12 K^2\mathbb P \big( \text{a given subsquare of area $\eta^2/K^2$ has 0 }
  \\ & \qquad \text{or more than $\lfloor L/(12 K^2) \rfloor$ points} \big). \numberthis\label{longunion}
\end{align*}
Let us fix $\eps>0$ and define 
\[ \gamma'(\lambda) = \frac{1}{12 M\lambda \eta^2(1+\eps)} \Big( \frac{\ell(\sqrt 5 \eta/K)}{\tau}-N_0\Big) .\numberthis\label{Omega1/lambdadef} \] Then for $\gamma=\gamma'(\lambda)$, we have $L=\lfloor 12 \lambda \eta^2(1+\eps) \rfloor$ in \eqref{Ndef}. Using this, \eqref{longunion}, and the stationarity of $X^\lambda$, it suffices to show that for all sufficiently large $\eta$, 
\[ \P \Big( 1 \leq \# \big(X^\lambda \cap Q_{\eta/K} \big) \leq  \lambda (1+\eps) \eta^2/K^2  \Big) \numberthis\label{goodregimeprobability} \]
tends to one as $\lambda \to \infty$. Indeed, then, using that $\gamma'(\lambda)$ defined in \eqref{Omega1/lambdadef} is $\Omega(1/\lambda)$, further that the set of edges of $g_{(\gamma)}(X^\lambda)$ is stochastically decreasing in $\gamma$, we can conclude the lemma. But for $\eta$ so large that $\Lambda(Q_{\eta/K})$ is also bounded away from zero,  the convergence of \eqref{goodregimeprobability} to zero  can be verified using an estimate analogous to \eqref{towerPalm}.
\ProofEnde

\section*{Acknowledgements}
The author thanks A.~Hinsen, C.~Hirsch, B.~Jahnel, W.~König, and R.~Löffler for interesting discussions and useful comments, and A.~Hinsen also for his help with creating Figure~\ref{figure-PVTgammas}.

\bibliographystyle{alea3}
\bibliography{alea_b2}

\end{document}